\documentstyle [12pt]{article}
\newcommand{\vsp}{\vspace{4 mm}}
\newcommand{\k}{\kappa}
\newcommand{\g}{g}
\newcommand{\A}{{\cal A}}

\newcommand{\D}{{\cal D}}
\newcommand{\Ho}{{\cal H}}
\newcommand{\h}{{\cal H}}
\newcommand{\Y}{{\cal Y}}
\newcommand{\ah}{a}
\newcommand{\Co}{\mbox{\boldmath $C$}}
\newcommand{\U}{{\cal U}}
\newcommand{\ra}{\rightarrow}
\newcommand{\Ra}{\Rightarrow}
\newcommand{\Lr}{\Leftrightarrow}
\newcommand{\be}{\begin{equation}}
\newcommand{\ee}{\end{equation}}
\newcommand{\nn}{\nonumber}
\newcommand{\bea}{\begin{eqnarray}}
\newcommand{\eea}{\end{eqnarray}}
\newcommand{\ot}{\otimes}
\newcommand{\op}{\oplus}
\newcommand{\C}{{\cal C}}
\newcommand{\Z}{\mbox{\boldmath $Z$}}
\newcommand{\va}{\mbox{\boldmath $1$}}
\newcommand{\vac}{\mbox{\boldmath $1$}}
\newcommand{\kj}{\left( \begin{array}{c} k \\ j \end{array} \right) }
\newcommand{\ri}{\left( \begin{array}{c} r \\ i \end{array} \right) }
\newcommand{\rri}{\left( \begin{array}{c} r \\ r-i \end{array} \right) }
\newcommand{\no}{\left( \begin{array}{c} n \\ 0 \end{array} \right) }
\newcommand{\mi}{\left( \begin{array}{c} m \\ i \end{array} \right) }
\newcommand{\io}{\iota}
\newcommand{\ep}{\epsilon}
\newcommand{\De}{\Delta}
\newcommand{\del}{\partial}
\newcommand{\un}{\underline}
\newcommand{\W}{W^{\infty /2}\tilde{\ell}}

\newcommand{\th}{\theta}
\newcommand{\na}{\nabla}
\newcommand{\we}{\wedge}
\newcommand{\br}{{[\; ,\; ]}_{SN}}
\newcommand{\dna}{D_{\nabla}}
\newcommand{\Om}{\Omega}
\newcommand{\Li}{{\cal L}}
\newcommand{\ad}{{\mbox{ad}}}
\newcommand{\td}{\tilde{\delta}}
\newcommand{\f}{\frac}

\begin{document}

\begin{center} {\bf ON SOME GENERALIZATIONS OF BATALIN-VILKOVISKY ALGEBRAS}
\\ {\bf F{\" u}sun Akman} \\ {\em Dept. of Mathematics, Cornell Univ., 
Ithaca, NY 14853-7901} \\ {\em  akmanf@math.cornell.edu}
\end{center}
\vsp

\vsp

\vsp

{\bf Abstract}
\vsp

{\it We define the concept of higher order differential operators on a general 
noncommutative, nonassociative superalgebra $\A$, and show that a vertex 
operator superalgebra (VOSA) has plenty of them, namely modes of vertex 
operators. A linear operator $\De$ is a differential operator of 
order~$\leq r$ if an inductively defined $(r+1)$-linear form 
$\Phi_{\De}^{r+1}$ with values in $\A$ is identically zero. These forms
resemble the multilinear string products of Zwiebach.
When $\A$ is a
``classical'' (i.e. supercommutative, associative) algebra, and $\De$ is an 
odd, square zero, second order differential operator on $\A$, $\Phi_{\De}^2$
defines a ``Batalin-Vilkovisky algebra'' structure on $\A$. Now that a second 
order differential operator makes sense, we generalize this notion to
any superalgebra with such an operator, and show that all properties of the
classical BV bracket but one continue to hold ``on the nose''.
As special cases, we provide several examples of classical BV
algebras, vertex operator BV algebras, and differential BV algebras. We also
point out connections to Leibniz algebras and the noncommutative homology
theory of Loday. Taking the generalization process one step further,
we remove all conditions on the odd operator $\De$ to examine the changes
in the basic properties of the bracket. We see that a topological chiral
algebra with a mild restriction yields a classical BV algebra in the
cohomology.
Finally, we investigate the quantum BV master equation for (i)~classical BV 
algebras, (ii)~vertex operator BV algebras, and (iii)~generalized BV algebras,
relating it to deformations of differential graded algebras.}
\vsp

\vsp

\vsp

{\bf 1 Introduction}
\vsp

The algebra of Batalin-Vilkovisky (BV) quantization {\bf [BV]} has been
studied from several viewpoints.
In its most abstract form, a BV algebra is a supercommutative associative
algebra $\A$ with an odd, square zero, second order differential operator
$\De$. From
these axioms follows the existence of an odd Lie bracket $\{ \; ,\; \}$ 
measuring the failure of $\De$ to be a derivation. The definition depends
very much on the notion of a second order differential operator, which seems 
to have been defined so far for operators on supercommutative, associative
algebras only. In particular J.-L. Koszul 
{\bf [Ko]} proposes a definition of higher order differential operators on
such an algebra $\A$ in terms of the vanishing of $r$-linear forms 
$\Phi_{\De}^r$ with values in $\A$, which is not applicable to other types 
of algebras. At this point we should make the meanings of ``super'' and
``odd'' very clear. A superalgebra $\A$ is simply an algebra which has
an integer grading preserved by the multiplication, and which is also 
expected to 
have identities differing from some ordinary version by certain powers of $-1$,
wherever the product of two homogeneous elements is written in two different
ways. For example, a supercommutative
(a.k.a. graded-commutative) algebra is a superalgebra where the identity
$ab-(-1)^{|a||b|}ba=0$ is satisfied for all homogeneous elements $a$, $b$.
An odd operator $\De$ on $\A$ is an operator which shifts the grading by an
odd integer $|\De |$ (mostly there is an implicit assumption of homogeneity).  
One may assume this odd integer is $\pm 1$ in a computation where it only
appears as a power (of $-1$), as is the case with most of this article.
$\De$ is called a superderivation (or simply a derivation) of $\A$ if it
satisfies the usual product rule with modifications; i.e. if
$\De(ab)=\De(a)b+(-1)^{|a||\De |}a\De(b)$ for every homogeneous $a$.
\vsp

Examples of BV algebras (which we may now call {\sl classical}) have been
floating around: The homology complex $(\wedge\g ,\del)$ of a finite 
dimensional Lie algebra $\g$ is one. Another geometric example involves the
Schouten-Nijenhuis bracket on the contravariant antisymmetric tensor
fields on a finite dimensional manifold ({\bf [Ko]}).
Lian and Zuckerman showed in {\bf [LZ1]} 
that the cohomology of a BRST complex with the Wick product is a BV algebra,
where the BV operator is the weight zero mode of the anti-ghost vertex
operator $b(z)$ (also see the follow-up by Penkava and Schwarz {\bf [PS]}).
In a related construction Bouwknegt and Pilch (in collaboration with
McCarthy) showed that a certain semi-infinite (BRST) cohomology of the
${\cal W}_3$ algebra has a BV structure, with the ``same'' BV operator
{\bf [BP]}. Most of these examples are subspaces or subquotients of richer
structures, some of which are vertex operator superalgebras (VOSA) or
close relatives.  Lian and Zuckerman study in {\bf [LZ1]} several properties
of the BV bracket which hold ``off-shell'', i.e. on the cochain complex
itself, which is a VOSA. One of our goals is to isolate -with a lot of 
hindsight- the most general properties of a BV bracket and elucidate
the algebra implicit in {\bf [LZ1]}.
\vsp

Some of the problems to be overcome in a broader definition of a BV
algebra are the following: The new definition must involve an
unambiguous notion of a second order differential operator for the most general
types of algebras. The corresponding BV bracket should retain most of 
its desirable properties. There should be abundant natural examples
among VOSA's, and every major existing example should be a subquotient
of one. And finally, these ``higher brackets'' had better have some
geometric or physical significance. Our treatment will fulfill most of
these criteria.
\vsp

In Section~2.1 we propose a definition of {\it higher order differential
operators}
on a superalgebra, which need not be supercommutative or associative.
A linear operator $\De$ on $\A$ is an $r$-th order differential operator if 
and only if a (recursively) well-defined $(r+1)$-linear form 
$\Phi_{\De}^{r+1}(a_1,\cdots ,a_{r+1})$ is identically zero. The
definition agrees with Koszul's for classical algebras. Section~2.2
has the important result (Theorem~2.2) that the modes $u_0$, $u_1$,
$u_2$, ... of a vertex operator $u(z)=\sum u_nz^{-n-1}$ are respectively
first, second, third, ... order differential operators on the VOSA with 
respect to the Wick product, so that any odd element $u$ of a VOSA with a 
square zero mode $u_1$ provides us with a BV operator!
\vsp

Section~3.1 contains the definition of classical Gerstenhaber 
and BV algebras, and in Section~3.2, we look at examples of classical
BV algebras (classical Lie algebra cohomology and homology complexes
with coefficients in a commutative associative algebra, more of the above
with multiplication and substitution operators as BV operators, the
Weil algebra, the cohomology of a topological chiral algebra, and skew
multivector fields with the Schouten-Nijenhuis bracket),
supplying simple proofs that the operators under consideration are 
second order differential operators. Section~3.2.1 provides an answer
to Exercise~10 in {\bf [LZ4]}. 
\vsp

We give the definition of a {\it generalized BV algebra (GBVA)} and its 
bracket (namely, an algebra $\A$ with an odd, square zero, second order
differential operator $\De$ and $\{ a,b\} =(-1)^{|a|}\Phi_{\De}^2(a,b)$) in 
Section~4.1.
A GBVA is possibly the most prominent nontrivial example of a ``Leibniz
algebra'', which is the main ingredient of the noncommutative homology 
theory of Loday {\bf [Lo]}.
(A Leibniz algebra is a generalized Lie algebra where the superderivation 
(Leibniz) property of the bracket holds, but skew-symmetry need not.)
We define a {\em vertex operator BV algebra (VOBVA)} to be a GBVA where 
$\A$ is a VOSA with the Wick product and $\De$ is some $u_1$.
Since many known examples of BV algebras can be obtained as
the cohomology of a generalized BV algebra, we define a {\it differential
BV algebra (DBVA)} as a GBVA which exhibits an additional odd, square
zero operator $D$ and a diagonalizable operator $L$ such that
$D\De +\De D=L$ ($L=0$ or $L=L_0$ of $Vir$ most of the time), with the
condition that the cohomology of $D$ is also a GBVA. This latter may
seem superfluous, but $D$ need not always be a derivation
and hence its cohomology need not be of the form a subalgebra modulo an 
ideal. We show that all properties of the BV bracket survive these
definitions, with the exception of skew-symmetry, which is modified
(Proposition~4.8). The extra term $(-1)^{|a|}\tilde{\Phi}_{\De}^2
(a,b)$ associated with the skew-symmetric product ${[a,b]}=ab-
(-1)^{|a||b|}ba$ vanishes if $\A$ is supercommutative. Next, we
remove all restrictions on $\De$ as well and observe that the
identities in Proposition~4.8 are modified by certain $\Phi^2$ and
$\Phi^3$ terms, all of which vanish when $\A$ is supercommutative,
but not necessarily associative, and $\De$ is square zero and of order
two (Proposition~4.13).
In Section~4.2 we explore two representatives of VOBVA's, one of which is
the {\it vertex operator Weil algebra (VOWA)} with a huge number
of BV operators. The second part is a short reminder of the
properties of the well-known BRST complex, or more generally, of a 
topological chiral algebra (TCA). Proposition~4.14 asserts that
a TCA where the field $G(z)$ (definition in 4.2.2) is primary 
collapses to a classical BV algebra as $G_0^2$ becomes zero in the
$Q$-cohomology.
\vsp

Section~5.1 deals with the {\em quantum BV master equation} for classical
BV algebras. We show why $\De(exp(\frac{i}{\hbar}W))=0$ follows from 
$\{ W,W\} =2i\hbar\De(W)$. In Section~5.2 we investigate the same equation 
for VOBVA's, and finally in~5.3, we discuss the meaning of the master
equation for GBVA's (keeping the conditions on $\De$ intact), as related 
to the deformation
theory of differential graded algebras. A solution of this equation is
interpreted as a deformation of the square zero, second order operator
$\De$ where all identities are preserved.
\vsp

We do not give the full definition, nor all the properties, of vertex operator
superalgebras, as there are very good (and increasingly accessible)
accounts in literature. Apart from the classics {\bf [FLM, FHL, DL]},
and comprehensive reviews like {\bf [Geb]}, we recommend the sequence
{\bf [LZ2-LZ4]} which stresses similarities with supercommutative
associative algebras and builds up the theory from scratch. 
It suffices at the moment to say that a VOSA
is a $\Z$-bigraded vector space $V$ (one $L_0$ and one super grading) where
each element $u$ is represented uniquely (and linearly) by a ``vertex
operator'' $u(z)=\sum u_nz^{-n-1}$, with $u_n\in End(V)$. There is an
action of the Virasoro algebra by $L(z)=\sum L_nz^{-n-2}$; the eigenvalues
of $L_0$ are bounded from below, and $L_{-1}$ acts as formal differentiation.
The vacuum element $\va$ is represented by $1\cdot z^0$ and is the identity
element with respect to the multiplication $u\times_{-1}v=u_{-1}v$ (the
{\em Wick product}). If the usual {\em Cauchy-Jacobi identity} involving the 
relations of modes $u_n$ becomes too oppressive,
one may visualize the alternative hidden in {\bf [DL]}, which says that for any
$u$, $v$ in the VOSA there exists a sufficiently large positive integer $t$ 
such that
\[ {[u(z_1),v(z_2)]}(z_1-z_2)^t=0 \]
as formal power series. We will frequently make use of the identities
(see e.g. {\bf [FLM]} or {\bf [Geb]})
\be (u_mv)_n=\sum_{i\geq 0}(-1)^i\mi (u_{m-i}v_{n+i}-(-1)^{m+|u||v|}
v_{m+n-i}u_i)\label{geb1}\ee
and    
\be {[u_m,v_n]}=\sum_{i\geq 0}\mi(u_iv)_{m+n-i}\label{geb2}\ee
for $m$, $n\in\Z$ ($\mi$ is given by $m(m-1)\cdots (m-i+1)/i!$
if $i\geq 1$, and $1$ if $i=0$). Note that a VOSA is {\sl not} thought of as a 
super-Virasoro module, as the name sometimes implies in the 
mathematical physics literature.
\vsp

We should also mention operadic constructions of BV algebras (such as
{\bf [Get1]}, {\bf [Hu]}). Recently (while this work was in progress)
Kimura, Stasheff, and Voronov proposed an abstract partial solution to
the problem of lifting the BV algebra structure in {\bf [LZ1]} to the 
cochain complex, by introducing ``filtered topological gravity'' whose
state space is a ``commutative homotopy algebra'', and stating that such
algebras can be obtained as algebras over a certain operad, with some
vanishing condition (Thm.~0.1 of {\bf [KSV]}). We will not study this
aspect of BV algebras at all, but point out that there is an obvious way
of writing a linear operad whose algebras are exactly the generalized
BV algebras as defined in our paper. Finding a topological operad whose
cohomology is this linear operad is another matter.
\vsp

\vsp

\vsp

{\bf 2 Higher Order Differential Operators}
\vsp

{\it 2.1 Definitions}
\vsp

There is a notion of ``higher order differential operators'' for 
(super) commutative, associative algebras, which is consistent with the 
idea of composites of partial derivatives on an algebra of functions.
Let  $\A =\op_j\A^j$ be  any superalgebra, not necessarily supercommutative
or associative. A linear map $\De :\A \ra \A$ is said to be {\em homogeneous
of (super)degree k}, written $|\De |=k$, if $\De :\A^j \ra \A^{j+k}$ for all 
$j\in \Z$. The map $\De$ is also said to be a {\em first order differential
operator} (or {\em derivation}) on $\A$ if
\be \De(ab)=\De(a)b+(-1)^{|a| |\De |}a\De(b). \label{prod} \ee
A {\em second order differential operator} is a map $\De$ such that the bracket
operator $\{ \;a, \; \}$ defined by
\be \{ a,b \}=\De(ab)-\De(a)b-(-1)^{|a||\De |}a\De(b) \label{sec} \ee
is a derivation for every homogeneous $a$. (The bracket measures the 
deviation of $\De$ from being a derivation.) J.-L.~Koszul defines in {\bf [Ko]}
linear maps (equivalently, $r$-linear forms)
\[ \Phi_{\De}^r : \A^{\ot r} \ra \A \; \; (r\geq 1) \]
for every linear operator $\De$ on a supercommutative associative 
algebra $\A$ by
\[ \Phi_{\De}^r (a_1,\cdots ,a_r)=m\circ (\De\ot id_{\A})\lambda^r(a_1 \ot
\cdots \ot a_r) \]
where $m(a\ot b)=ab$ is the multiplication map, and
\[ \lambda^r(a_1\ot\cdots\ot a_r)=(a_1\ot 1-1\ot a_1)\cdots (a_r\ot 1-1\ot
a_r), \]
the right hand side being a product in the supercommutative algebra $\A\ot\A$.
Then
\bea \Phi_{\De}^1 (a)=&&\De(a)-\De(1)a \nn \\ \Phi_{\De}^2 (a,b)
=&&\De(ab)-\De(a)b-(-1)^{|a||b|}
\De(b)a+\De(1)ab \label{ko} \\ \Phi_{\De}^3 (a,b,c)=&&\De(abc)-\De(ab)c
-(-1)^{|a|(|b|+|c|)}\De(bc)a \nn \\ && -(-1)^{|c|(|a|+|b|)}\De(ca)b
+\De(a)bc+(-1)^{|a|(|b|+|c|)}\De(b)ca \nn \\ && +(-1)^{|c|(|a|+|b|)}
\De(c)ab-\De(1)abc \nn \\ && \vdots \nn \eea
and a map $\De$ (of degree $k$)is said to be a {\em differential operator on} 
$\A$ {\em of order}~$\leq r$ (written $\De \in {\D}_{r}^k$) iff 
$\Phi_{\De}^{r+1}$ is identically zero. The subspaces $\D_r^j$ of $Lin(\A)$ 
satisfy
\bea && \mbox{(i)} \; \; \D_1^j\subset\D_2^j\subset\cdots\subset\D_r^j
\subset\D_{r+1}^j
\subset\cdots \nn \\ && \mbox{(ii)} \; \; \D_r^j\D_s^k\subset\D_{r+s}^{j+k}
\label{bir} \\ && \mbox{(iii)} \; \; 
{[\D_r^j,\D_s^k]}\subset\D_{r+s-1}^{j+k}. \nn \eea

{\bf Remark 2.1} {\em We will consider only differential operators which are
``graded'', i.e. which are finite linear combinations of elements of the} 
$\D_r^j$.
\vsp

Koszul asserts that every $\Phi_{\De}^{r+1}$ can be written in terms of
$\Phi_{\De}^r$, such as
\[ \Phi_{\De}^3(a,b,c)=\Phi_{\De}^2(a,bc)-\Phi_{\De}^2(a,b)c-(-1)^{|b||c|}
\Phi_{\De}^2(a,c)b, \]
so that property (i) above is trivial (no other inductive formula is given).
\vsp

In order to generalize the notion to a noncommutative, nonassociative
superalgebra $\A$ 
(e.g. a VOSA) this last formula must be studied. A good inductive definition
of $\Phi_{\De}^r$ should involve only one multiplication at a time,
and no arbitrary change of order of arguments. We propose
a new recursive definition:\footnote{Prof. Koszul commented in a private
communication that there is an older recursive definition due to
Grothendieck {\bf [Gr]}: a linear operator $\De :\A\rightarrow\A$
is of order $\leq r$ if and only if the (super)commutator $[\De ,m_a]$ is of
order $\leq r-1$ for all $a\in\A$, 
where $m_a$ is left multiplication by $a$, and operators
of order $\leq -1$ are zero. This definition is different from mine and does
not work for VOSA's.}
\bea \Phi_{\De}^1 (a)=&&\De(a) \nn \\ \Phi_{\De}^2 (a,b)=&&\Phi_{\De}^1(ab)-
\Phi_{\De}^1(a)b-(-1)^{|a||\De |}a\Phi_{\De}^1(b) \nn \\  
\Phi_{\De}^3 (a,b,c)=&&\Phi_{\De}^2(a,bc)
-\Phi_{\De}^2(a,b)c-(-1)^{|b|(|a|+|\De |)}b\Phi_{\De}^2(a,c) 
\label{newphi} \\ && \vdots \nn \\ 
\Phi_{\De}^{r+1}(a_1,\cdots ,a_{r+1})=&&\Phi_{\De}^r(a_1,\cdots ,a_r a_{r+1})
-\Phi_{\De}^r(a_1,\cdots ,a_r)a_{r+1} \nn \\ && -(-1)^{|a_r|(|a_1|
+\cdots +|a_{r-1}|+|\De |)} a_r\Phi_{\De}^r(a_1,\cdots ,
a_{r-1},a_{r+1}) \nn \\ && \vdots \nn \eea
which fits the bill. The extra powers of $(-1)$ come from 
the rearrangement of the symbols $\De$, $a_1$, ..., $a_{r+1}$.
Note that each $\Phi_{\De}^{r+1}(a_1,\cdots ,a_{r+1})$ gives the
deviation of $\Phi_{\De}^r(a_1,\cdots ,a_{r-1},\; \; \; )$ (which we can
interpret as a {\em higher bracket}; see {\bf [LS, Zw, KSV]} on ``strongly
homotopy Lie algebras'', ``multilinear string products'', and related notions) 
from being a derivation, where the order
in the definition of the product rule is fixed by~(\ref{prod}).
Then carefully keeping track of the order of multiplication and of the 
arguments $a_i$, we see that
\bea \Phi_{\De}^1(a)=&&\De(a) \nn \\ \Phi_{\De}^2(a,b)=&&\De(ab)-\De(a)b-
(-1)^{|a||\De |}a\De(b)
\label{newsec} \\ \Phi_{\De}^3(a,b,c)=&&\De(a(bc))-(-1)^{|a||\De |}a\De(bc)
\nn \\ && -(-1)^{|b|(|a|+|\De |)}b\De(ac)-\De(ab)c
-\De(a)(bc)+(\De(a)b)c \nn \\ && +(-1)^{|b|(|a|+|\De |)}b(\De(a)c)
+(-1)^{|a||\De |}(a\De(b))c+(-1)^{|a||b|+|\De |(|a|+|b|)}
b(a\De(c)) \nn \\ && \vdots \nn \eea
leads to the correct (or plausible, unambiguous, ...) definition of higher
order differential operators in a general algebra $\A$. If $\A$ happens to
have an identity (such as a VOSA and its vacuum element, which is the identity
with respect to the Wick product), one may wish to replace the first
line in (\ref{newphi}) by
\be \Phi_{\De}^1(a)=\De(a)-\De(1)a \label{mod} \ee
and rewrite (\ref{newsec}) accordingly.
In the polynomial algebra $\Co [x]$, for example, multiplication by a 
polynomial $p(x)$ would be a differential operator of order zero, 
and would annihilate all forms. The first order differential operators would
be of the form $p(x)\frac{d}{dx}$ and annihilate all except possibly
the first form, etc.
\vsp

Definitions (\ref{ko}) and (\ref{newphi})-(\ref{mod}) have been checked to
agree up to and including $\Phi_{\De}^4$ for classical algebras, and it
should not be too difficult to give a general proof of equivalence in this
special case. The following familiar formula for $\Phi_{\De}^4$ clearly
shows the best way to write any expression $\Phi_{\De}^r$ in a classical
algebra. We apply the operator $\De$ to $k$ entries out of $r$ 
($1\leq k\leq r$) in every possible way (up to order) and multiply each of
these
expressions with the remaining $r-k$ entries. We then add up all terms,
modifying each by $(-1)^{r-k}$ {\sl and} by the sign of the permutation
of the symbols $\De$, $a$, $b$, $c$, ... under consideration:
\bea \Phi_{\De}^4(a,b,c,d)=&&\De(abcd)\nn \\
&& -(-1)^{|\De ||a|}a\De(bcd)-(-1)^{|b|(|a|+|\De |)}b\De(acd)\nn \\
&& -(-1)^{|c|(|a|+|b|+|\De |)}c\De(abd)-\De(abc)d \nn \\
&& +(-1)^{|\De |(|a|+|b|)}ab\De(cd)+(-1)^{|\De ||a|}a\De(bc)d\nn \\
&& +(-1)^{|b|(|a|+|\De |)}b\De(ac)d+(-1)^{|c|(|\De |+|b|)+|\De ||a|}ac\De(bd)
\nn \\
&& +(-1)^{(|b|+|c|)(|\De |+|a|)}bc\De(ad)+\De(ab)cd\nn \\
&& -\De(a)bcd-(-1)^{|\De ||a|}a\De(b)cd\nn \\
&& -(-1)^{|\De |(|a|+|b|)}ab\De(c)d-(-1)^{|\De |(|a|+|b|+|c|)}abc\De(d).
\nn\eea

\vsp

\vsp

{\it 2.2 Higher Order Differential Operators on VOSA's}
\vsp

We will show that in a vertex operator 
superalgebra, the modes $u_0$, $u_1$, $u_2$,~... of any vertex operator 
\be u(z)=\sum_{n\in \Z} u_nz^{-n-1}, \label{uz} \ee
namely the coefficients of $z^{-1}$, $z^{-2}$, $z^{-3}$,~..., are
differential operators
on the VOSA of orders~$\leq$ 1, 2, 3, ... respectively (they annihilate
the vacuum). Meanwhile all the other coefficients are ``left
multiplications'', i.e. they annihilate the modified expression~(\ref{mod}).
\vsp

For vertex operators written in the standard form~(\ref{uz}), the 
infinitely many multiplications
\be a\times_n b=a_n b=a_n b_{-1}\va\label{m} \ee
are given by
\[ a\times_n b=Res_z z^{n}a(z)b, \]
and $\times_{-1}$ is the ``normal ordered product'', a.k.a. the ``Wick
product''. It is well known that a residue, $u_0$, is a derivation of the
VOSA with respect to all products $\times_n$:
\[ u_0(a_n b)=(u_0 a)_n b+(-1)^{|u||a|}a_n (u_0 b) \]
(compute ${[u_0,a_n]}$ from (\ref{geb2}) and apply to $b$). Then
\[ \Phi_{u_0}^2(a,b)\equiv 0 \]
with respect to all $\times_n$, $n\in \Z$, that is, $u_0 \in \D_1$. However
at this point we restrict the definition of a differential operator to
the case of the Wick product, since the other modes are not as well-behaved 
as the residue. Let $ab$ denote $a\times_{-1} b$ from now on.
\vsp

{\bf Theorem 2.2} {\em In a VOSA $V$, the modes $u_n$ of a vertex 
operator~(\ref{uz}) for $n\geq 0$ are higher order differential operators on
$V$, namely, 
\[ u_n\in \D_{n+1} \; \; \;\mbox{for}\; \; \; n\geq 0, \]
or equivalently,
\[ \Phi_{u_n}^{n+2}\equiv 0 \;\; \; \mbox{for}\; \; \; n\geq 0. \]
The remaining modes are left multiplications, that is,}
\[ \Phi_{u_n}^k \equiv 0 \; \;\; \forall n\leq -1, \;\; k\geq 1. \]
\vsp

{\bf Remark 2.3} {\em Note that condition (i) of (\ref{bir}) is automatically
satisfied.}
\vsp

{\em Proof of Theorem 2.2.} For $r\geq 1$, assume $u_0, \cdots , u_{r-1}$ 
have been shown to be
in $\D_1 , \cdots ,\D_r$ respectively for all $u$. We start with
\[ \Phi_{u_r}^1(a) = u_r a \; \; \; \; \mbox{(definition)}, \]
so that
\bea \Phi_{u_r}^2(a,b) &&= u_r(a_{-1}b)-(u_r a)_{-1}b-(-1)^{|u||a|}a_{-1}
(u_r b) \nn \\ &&=-(u_ra)_{-1}b+[u_r,a_{-1}]b\nn \\
&&= -(u_ra)_{-1}b+\sum_{i=0}^r\ri (u_ia)_{r-1-i}b \;\;\;
\mbox{(from (\ref{geb2}))} \nn \\
&&=-(u_ra)_{-1}b+\sum_{i=0}^r \rri (u_ia)_{(r-i)-1}b \nn \\
&&=-(u_ra)_{-1}b+\sum_{i=0}^r \ri (u_{r-i}a)_{i-1}b \;\;\;
\mbox{(replacing $i$ by $r-i$)}\nn \\
&&=\sum_{i=1}^r\ri (u_{r-i}a)_{i-1}b\label{del} \eea
where the subscripts of the operators applied to $b$ range from $0$ to
$r-1$, and by the induction step they lie in $\D_1,\cdots ,\D_{r}$, or
simply in $\D_r$. Hence
\[ \Phi_{u_r}^2(a,b)=\De(b)=\Phi_{\De}^1(b) \]
for some $\De\in\D_r$, keeping $a$ fixed (read $\De$ from (\ref{del})).
Now by an independent and 
easy induction we can show that
\[ \Phi_{u_r}^k(a_1,\cdots ,a_k)=\Phi_{\De}^{k-1}(a_2,\cdots ,a_k) \; \;
\mbox{for} \; \; k\geq 2 \]
and in particular
\[ \Phi_{u_r}^{r+2}(a_1,\cdots ,a_{r+2})=\Phi_{\De}^{r+1}(a_2,\cdots ,a_{r+2})
\equiv 0, \]
($a_1=a$, fixed) as $\De\in \D_r$.
Finally, the statement on left multiplications follows from
\bea \Phi_{u_n}^1(a) &&=u_na-(u_n\va)_{-1}a\nn \\
&&=u_na-\no u_na=0 \;\;\; \mbox{(from (\ref{geb1}))}\nn \eea
and the last Remark. $\Box$
\vsp

{\bf Remark 2.4} {\it Compositions and brackets of modes satisfy properties
similar to (ii) and (iii) in (\ref{bir}).}
\vsp

\vsp

\vsp

{\bf 3 Classical Batalin-Vilkovisky Algebras}
\vsp

{\it 3.1 Definitions and Properties}
\vsp

A {\em Gerstenhaber algebra} ({\bf [Ger]}) is a supercommutative, 
associative algebra
\[ \A = \op_j {\A}^j \]
equipped with an odd bracket $\{ \; \; , \; \; \}$ satisfying
\[ \{ a,b\} =-(-1)^{(|a|+1)(|b|+1)}\{ b,a \} \]
(skew-symmetry in the associated, shifted-graded super Lie algebra
$\hat{\A}=\op_j {\A}^{j+1}$), and
\[ \{ a, \{ b,c \} \} = \{ \{ a,b \} ,c \}+(-1)^{(|a|+1)
(|b|+1)} \{ b, \{ a,c \} \} \]
(the Leibniz rule, or superderivation property, in $\hat{\A}$), as
well as the superderivation (Poisson) rule with respect to the multiplication 
in $\A$:
\be \{ a,bc \} =\{ a,b \} c+(-1)^{(|a|+1)|b|}b\{ a,c\}. 
\ee
This last condition implies that bracketing with a homogeneous element
$a$ is a superderivation on $\A$ which changes the grading by $|a|$ plus an 
odd integer.
\vsp

A {\em Batalin-Vilkovisky algebra}, or {\em BV algebra}, is a Gerstenhaber
algebra where the bracket $\{ \; , \; \}$ is obtained from an odd, 
square zero, second order differential operator $\De$ (usually of superdegree
$\pm 1$):
\bea \{ a,b \} && =(-1)^{|a|}\Phi_{\De}^2(a,b)\label{above} \\ && =(-1)^{|a|}
\De(ab)-(-1)^{|a|}\De(a)b-a\De(b) \nn \eea
({\bf [BV]}, {\bf [LZ1]}, {\bf [Ko]}, {\bf [PS]}, {\bf [Get1]}). We will show 
in Section~4.1 that the identities above are modified by some $\Phi^2$ and
$\Phi^3$ terms when $\{ a,b\}$ is defined via (\ref{above}) for any 
superalgebra $\A$ and an odd operator $\De$ with no restrictions.
\vsp

\vsp

\vsp

{\it 3.2 Classical Examples}
\vsp

We will look at
\vsp

(i) $(\ah \ot \wedge \g',d)$: The classical Lie algebra cohomology complex
for a finite dimensional Lie algebra $\g$, and 
a commutative algebra $\ah$ on which $\g$ acts by derivations;

(ii) $(\ah \ot \wedge \g ,\del)$: The classical Lie algebra homology
complex for $\g$, $\ah$ as above;

(iii) Multiplication and substitution operators in (i) and (ii);

(iv) $S\g'\ot\wedge\g'$: The Weil algebra on $\g$ with several differentials;

(v) The cohomology of a BRST complex , or more generally,
the cohomology of a topological chiral algebra;

(vi) The skew multivector fields $A(M)$ on a finite dimensional paracompact
differentiable manifold $M$ with the Schouten-Nijenhuis bracket;
\vsp

and give simple arguments as to why the given operators are second
order differential operators. See also {\bf [Ko]}, {\bf [AKSZ]}, {\bf [Wi1]}, 
and {\bf [Sc]}.
\vsp

{\bf Remark 3.1} {\em Although most properties apply to a general Lie algebra,
we will restrict ourselves to finite dimensional semisimple $\g$ in this
section.}
\vsp

\vsp

\vsp

{\it 3.2.1 Classical Lie Algebra Cohomology and Homology}
\vsp

There is a unified ``semi-infinite'' construction of finite dimensional
Lie algebra cohomology and homology complexes ({\bf [A1, A2]}) which is
summarized below. We think of both complexes (i) and (ii) above as modules
over the associative superalgebra
\[ \Y\g =\U g \ot \C \g \]
where $\U\g$ is the universal enveloping algebra of $\g$, and $\C\g$ is the
Clifford algebra on $\g\op\g'$, the symmetric bilinear product being given
by the pairing between a fixed homogeneous basis $\{ \io(e_i)\}$ of $\g$ 
and its (restricted, i.e. piecewise) dual
$\{ \ep(e_i')\}$ in $\g'$. The superdegrees of the generators of types 
$\io$ and $\ep$ are taken to be $-1$ and $+1$ respectively. Every 
(co)homology related operator is realized as an inner derivation of $\Y\g$,
which then acts on (i) or (ii). Note that $\wedge\g'$ is spanned by wedge
products of $\ep(e_i')$'s and $\wedge\g$ by wedge products of $\io(e_i)$'s. 
There is an overall $\g$ action given by the operators $\th(x)\in Der(\Y\g)$,
namely
\[ \un{\th}(x)=\pi(x)+\un{\rho}(x)=\pi(x)+\sum_i:\io([x,e_i])\ep(e_i'):
\in \Y\g , \]
where $x\in\g$, $\pi(x)\in\U\g$, $\un{\rho}(x)\in\C\g$, and the notation
$:\;\; :$ (``normal ordered product'') simply means
\[ \un{\rho}(x)=\left\{ \begin{array}{ll}
-\sum_i\ep(e_i')\io([x,e_i]) & \mbox{in case (i)} \\
\sum_i\io([x,e_i])\ep(e_i') & \mbox{in case (ii).} \end{array} \right. \]
We say the $\io$'s are ``annihilation operators'' whereas the $\ep$'s are 
``creation operators'' for the cohomology complex (and vice versa for the
homology complex).  
It was shown in {\bf [A1]} that there exists a unique derivation $D$ 
($=d$ or $\del$) of $\Y\g$ of superdegree $+1$, satisfying
\[ D\; \io(x)=\un{\th}(x) \; \; \; \forall x\in\g. \]
Then $D$ turns out to be an inner derivation, with formula
\bea \un{D} &&= \sum_i\pi(e_i)\ep(e_i')+\sum_{i<j}:\io([e_i,e_j])
\ep(e_j')\ep(e_i'): \; \; \; \in\Y\g \nn \\
&&= \left\{ \begin{array}{ll} 
\un{d}=\sum_i\pi(e_i)\ep(e_i')+\sum_{i<j}\ep(e_j')\ep(e_i')\io([e_i,e_j])
& \mbox{in case (i)} \\ \un{\del}=\sum_i\pi(e_i)\ep(e_i')+\sum_{i<j}
\io([e_i,e_j])\ep(e_j')\ep(e_i') & \mbox{in case (ii).} \end{array}
\right. \label{dee}\eea
The characterization of $D$ above translates into the famous {\em Cartan
identity}
\[ \un{D}\io(x)+\io(x)\un{D}=\un{\th}(x) \]
relating operators acting on (i) or (ii). This identity favors $\io(x)$
over $\ep(x')$. We have the additional identities
\[ \un{D}^2=0 \]
and
\[ {[ \un{D},\un{\th}(x)]}=0 \; \; \; \forall x\in\g.\]
For both (i) and (ii) (taking $\ah =\Co$, the trivial representation,
for the moment, so that $\pi(e_i)=0$), the 
operator $\un{\rho}(x)$ is also a derivation
of the associative algebra $\wedge\g'$ or $\wedge\g$ for all $x\in\g$, as
\vsp

1)${[\un{\rho}(x),\; \; ]}$ is a derivation of $\C\g$,

2) $\un{\rho}(x)$ sends creation operators to creation operators,

3) $\un{\rho}(x)\vac =0.$
\vsp

Of course, from a simpler viewpoint, $\un{\rho}$ is nothing but the
extension of the (co)adjoint action of $\g$ to the exterior algebra
by derivations.
\vsp

How do we get a BV algebra? In case (i), $\un{d}$ is a derivation of
$\wedge\g'$ for similar reasons:
\vsp

1) ${[\un{d},\; \;]}$ is a derivation of $\C\g$,

2) $\un{d}$ sends $\ep$'s to expressions of type  $\sum\ep\ep$ (see
{\bf [A1]}),

3) $\un{d}\vac =0$.
\vsp

Also $\un{d}^2=0$, so $\un{d}$ is a BV operator on $\wedge\g'$ with
trivial BV bracket $\{ \; \; ,\; \; \}$.
\vsp

In case (ii), $\un{\del}$ is {\sl not} a derivation of $\wedge\g$, for it fails
to send creation operators to (products of) creation operators in $\C\g$:
\be \del\; \io(x)=\un{\rho}(x). \ee
However, $\un{\del}$ is a sum of genuine second order differential operators
on $\wedge\g$:
\[ \un{\del}=\frac{1}{2}\sum_i\un{\rho}(e_i)\ep(e_i') \]
(both $\un{\rho}(e_i)$ and $\ep(e_i')$ are derivations!), hence
the square zero operator $\un{\del}$ becomes a nontrivial BV operator.
\vsp

The generalization to a nontrivial module $\ah$ is easy. Case (i): The
second term of $\un{d}$ in (\ref{dee}) is already a derivation of $\wedge\g'$,
and the first satisfies
\bea && \sum_i\pi(e_i)\ep(e_i')(ab)=\sum_i\ep(e_i')\pi(e_i)(ab)\nn \\
&&= \sum_i\ep(e_i')((\pi(e_i)a)+a(\pi(e_i)b)) \nn \\
&&= ((\sum_i\ep(e_i')\pi(e_i))a)b+a((\sum_i\ep(e_i')\pi(e_i))b), \nn \eea
so that $\un{d}$ is a derivation.
\vsp

Case (ii): Once again, $\un{\del}$ is not a derivation in general, but is
a sum of second order differential operators:
\[ \un{\del} =\sum_i(\pi(e_i)+\frac{1}{2}\un{\rho}(e_i))\ep(e_i'). \]
\vsp

\vsp

\vsp

{\it 3.2.2 Multiplication and Substitution Operators in (i) and (ii)} 
\vsp

The odd operators $\io(x)$ and $\ep(x')$ which satisfy the Clifford algebra
relations
\bea && \io(x)\io(y)+\io(y)\io(x)=0 \nn \\
&& \ep(x')\ep(y')+\ep(y')\ep(x')=0\nn \\
&& \ep(x')\io(y)+\io(y)\ep(x')=x'(y)\cdot 1 \nn \eea
give rise to -often extremely trivial- BV algebra structures on the complex
and/or its (co)homology. Both are square zero for starters. Let us look at
$\io(x)$ first.
\vsp

In case (i), ${[\io(x), \; \; ]}$ is a derivation of $\C\g$; $\io(x)\vac
=0$; and ${[\io(x),\ep(y')]}\in\Co$; hence $\io(x)$ is a derivation of
$\wedge\g'$ and a (fake) BV operator.
\vsp

In case (ii), all desired properties hold, including sending the creation
operators ($\io$'s) to oblivion, but now $\io(x)\vac\neq 0$. That's where
we have to pass to the cohomology to get one of these trivial (zero) BV
operators, provided that the homology $H_{\ast}(\g,\Co)$ does have a 
supercommutative, associative algebra structure. Since $\g$ is semisimple,
this condition is satisfied:
\[ H_{\ast}(\g,\Co)=(\wedge\g)^{\g}, \]
and also
\[ H_1(\g,\Co)=\g /[\g ,\g]=0, \]
so that 
\[ \io(x)\vac =0 \; \; \; \forall x\in\g \; \; \; \mbox{in} \; \; \;
H_{\ast}(\g,\Co). \]
Then
\[ \io(x)\equiv 0 \; \; \; \mbox{on} \; \; \; H_{\ast}(\g,\Co). \]

{\bf Remark 3.2} {\em The anti-ghost operator $b_0$ in a BRST complex 
satisfies 
the Cartan identity like $\io(x)$, but it combines characteristics of
$\un{\del}$ and $\io(x)$.}
\vsp

Now, about $\ep$'s. Case (i): ${[\ep(x'),\; \; ]}$ is a derivation of
$\C\g$, ${[\ep(x'),\ep(y')]}=0$, $\ep(x')^2=0$, but $\ep(x')\vac\neq 0$.
Again for finite dimensional semisimple (even reductive) $\g$, we have
\[ H^{\ast}(\g,\Co)=(\wedge\g')^{\g}, \]
an algebra,
\[ H^1(\g,\Co)=0, \]
and
\[ \ep(x')\vac =0. \]
Then
\[ \ep(x')\equiv 0 \; \; \; \mbox{on} \; \; \; H^{\ast}(\g,\Co). \]
Case (ii): $\ep(x')$ is similarly shown to be a derivation
on $\wedge\g$, with $\ep(x')^2=0$.
\vsp

\vsp

\vsp

{\it 3.2.3 The Weil Algebra}
\vsp

The {\it classical Weil algebra} 
\[ W\ell =S\ell'\ot\wedge\ell'\;\;\;
\mbox{($\ell=$ a finite dimensional Lie algebra)}\]
is a well-known object in differential geometry, used for example as a 
model for connections and curvature, and to define equivariant forms
{\bf [MQ]}. A deep exposition of its algebraic properties can be found in
{\bf [GHV]}. The Weil algebra with one of the differentials, $\un{d}$, 
is also an example of a 
classical Lie algebra cohomology complex (type~(i) in Section~3.2). We will
only point out that we obtain trivial BV brackets on the (cohomology) BV
algebras
\[ (W\ell ,\un{d}),\;\;\; (W\ell ,\un{h}),\;\;\; (W\ell ,\un{\k}),\;\;\;
(W\ell ,\un{d}+\un{h})\]
where the four differentials (resp. the Lie algebra cohomology differential,
the Koszul differential, a homotopy operator, and the Weil differential)
are all first order differential operators (see {\bf [A2]}). 
However, a ``change of vacuum'' results in the homology complex (case~(ii)
of~3.2)
\[ (W\ell)'=S\ell\ot\wedge\ell \]
with a reversal of multiplication and substitution operators, and 
normal ordering changes accordingly (annihilation operators
to the right of creation operators, as usual). The $\un{\del}$ homology is
\[ H_{\ast}((W\ell)',\un{\del})=(S\ell)^{\ell}\ot(\wedge\ell)^{\ell}
\;\;\;\mbox{for reductive $\ell$}, \]
where $\un{\del}$ is now a second order differential operator giving rise to a
nontrivial BV bracket in general. 
\vsp

\vsp

\vsp

{\it 3.2.4 The Cohomology of a Topological Chiral Algebra}
\vsp

The reader may want to come back to this section after reading the vertex 
operator
version in Section~4.2.2. We will assume the definitions and properties of
topological chiral algebras and BRST complexes in~4.2.2 for now, as well as
\vsp

{\bf Lemma 3.3} {\bf ([Li])} {\em Let
\[ V=\op_{n\geq 0}V[n] \]
be a VO(S)A with nonnegative weights only. Then $V[0]$ is a 
(super)commutative, associative algebra under the Wick product.}
\vsp

On a BRST complex $V$ with differential $Q$, we have the ``semi-infinite
Cartan identity''
\be Qb_0+b_0Q=L_0 , \label{sici}\ee
where $b_0$, $L_0$ are not residues, but coefficients of $z^{-2}$ of the
vertex operators
\[ b(z)=\sum b_nz^{-n-2} \; \; \; \mbox{(anti-ghost)} \]
and
\[ L(z)=\sum L_nz^{-n-2} \; \; \; \mbox{(stress-energy tensor; Virasoro),}\]
written in a non-standard way when compared with (\ref{uz}). The operator
$b_0$ plays the role of $\io(L_0)$ and is odd, square zero, and most 
importantly a second order differential operator on $V$, as was shown. The
cohomology
of this complex, namely $\Ho^{\ast}=Ker(Q)/Im(Q)$, was first shown by Lian
and Zuckerman to be a classical BV algebra with $\De =b_0$. We provide here
a short proof. 
\vsp

{\bf Proposition 3.4 ([LZ1])} {\em The cohomology $\h^{\ast}$ of a BRST
complex is a classical BV algebra under the induced Wick product and the
operator $b_0$.}
\vsp

{\em Proof.} $\h^{\ast}$ is a VOSA as $Q$ is a residue and $L(z)\in
Ker(Q)$. Thanks to (\ref{sici}), $\h^{\ast}$ is also confined to weight
zero. But then $\h^{\ast}$ is a VOSA with trivial $L_0$ grading, hence by
Lemma~3.3 it is a supercommutative associative algebra with respect to the
Wick product. The operator $b_0$, a second order differential operator on the 
VOSA according to Theorem~2.2, acts on $\h^{\ast}$ as explained later in 
Proposition~4.7. $\Box$
\vsp

{\bf Remark 3.5} {\em The commutative subalgebra $\h^0$ is what Witten
calls the {\em ground ring} of the theory {\bf [Wi2]}.}
\vsp

The cohomology of a topological chiral algebra is similarly a classical
BV algebra provided that the operator $G_0$ (analogue of $b_0$) is
square zero, at least in the cohomology. We will discuss conditions
on $G(z)$ in~4.2.2.
\vsp

We will see that at least algebraically it is more natural to consider the 
whole BRST complex as a ``vertex operator BV algebra'',and to obtain
$(\Ho^{\ast},b_0)$ as a classical subquotient of $(V,b_0)$. Note that
in Section~4.2.2 we denote a BRST complex by $V\ot \we^{\ast}$ and not
by $V$. In both
{\bf [LZ1]} and {\bf [PS]} there are indications that ``higher brackets''
on the original complex would be interesting.
\vsp

\vsp

\vsp

{\it 3.2.5 The Schouten-Nijenhuis Bracket on Skew Multivector Fields}
\vsp

The Schouten-Nijenhuis bracket $\br$ ({\bf [Sn, Ni, Ko, Mi2]}) is a 
differential invariant (concomitant) of skew multivector
fields which generalizes the Lie bracket of two vector fields. The bracket
satisfies skew-symmetry as well as the Jacobi and Poisson identities (with
respect to the properly shifted gradings) and makes
\[ A(M)=\mbox{Sec}(\we TM)=\we_{C^{\infty}(M)}\mbox{Sec}(TM) \]
into a Gerstenhaber algebra ($M$ is an $n$-dimensional paracompact
differentiable manifold). The definition is given by
\bea && \br :A^p(M)\times A^q(M)\ra A^{p+q-1}(M)\nn \\
&& {[X_1\we\cdots\we X_p,Y_1\we\cdots\we Y_q]}_{SN} \nn \\
&& =\sum_{i=1}^n\sum_{j=1}^n(-1)^{i+j}[X_i,Y_j]\we X_1\we\cdots\we
\hat{X_i}\we\cdots\we X_p\we Y_1\we\cdots\we\hat{Y_j}\we\cdots\we Y_q
\;\;\;\;\mbox{for $p$, $q\geq 1$} \nn \eea
and
\[ {[f,u]}_{SN}=\io(df)u\;\;\;\;\;\mbox{for $f$}\in A^0(M)=C^{\infty}(M),
\;\;\;\; u\in A(M). \]
Koszul shows that $(A(M),\br )$ can moreover be made into a BV algebra:
\vsp

{\bf Proposition 3.6 ([Ko])} {\it For every torsion free linear connection
$\na$ in $TM$, there exists a unique differential operator $\dna$ which
generates $\br$ via $\Phi_{\dna}^2$, such that}
\[ \dna(X)=-\mbox{div}_{\na}(X)=-\mbox{Tr}(Y\mapsto [X,Y]-\na_X(Y)).\]
{\it This map is the composition of two maps: One sending $\na$ to the
induced connection in $\we^nTM$, and the other the canonical bijection
of the affine space of linear connections in $\we^nTM$ onto the set of
differential operators which generate the Schouten-Nijenhuis bracket. In
this bijection the images of connections of zero curvature in $\we^nTM$
are square zero operators.}
\vsp

We will omit the details of the proof; however, we will give a local
expression for $\dna$ (on an open set $U\subset M$):
\[ D_U=-\sum_{i=1}^n\io(\alpha_i)\na_{X_i},\]
where $X_1$,...,$X_n$ is a basis for the module $A^1(U)$ of vector
fields on $U$, $\alpha_1$,...,$\alpha_n$ is the dual basis of 1-forms,
and $\io(\alpha_i)$ are ``substitution operators'', or ``interior products''
(analogues of the substitution operators $\io(X)$ acting on differential
forms, i.e. on $\mbox{Sec}(\we T^{\ast}M)$).
We obtain a differential operator $\dna$ of degree~$-1$ and order~$\leq 2$
by patching. There are two obvious examples where $\dna$ is square zero.
One occurs when $\na$ is the Levi-Civita connection associated to a 
Riemannian metric on $M$. The second one is the case of an orientable
manifold $M$ with a nonvanishing section $v$ of $\we^nTM$, which then
determines a flat connection in $\we^nTM$ whose image is a square zero
operator generating the Schouten-Nijenhuis bracket.
\vsp

It would be of interest to know which bilinear differential concomitants
(such as the brackets studied by Fr\"{o}licher, Nijenhuis, Richardson,
and Schouten, and the generalized Poisson bracket mentioned at the end of
Section~4.1) are BV brackets, even though the differential forms under
consideration may have other than a supercommutative, associative algebra
structure. We recommend {\bf [KMS]} for background and references.
\vsp

\vsp

\vsp

{\bf 4 Generalized BV Algebras}
\vsp

{\it 4.1 Definitions and Properties}
\vsp

The existence of higher order differential operators and higher brackets on a
VOSA $V$ leads naturally to the concept of a {\em vertex operator 
Batalin-Vilkovisky algebra}. Even more generally, 
\vsp

{\bf Definition 4.1} A {\em generalized Batalin-Vilkovisky algebra (GBVA)}
is a superalgebra $\A$ (possibly with identity, $1_{\A}$) with an odd,
square zero, second order differential operator $\De$ (which annihilates the
identity, if any). The {\em generalized BV bracket} $\{ \; , \; \}$ is given by
\be \{ a,b \} =(-1)^{|a|}\Phi_{\De}^2(a,b)\label{bracket} \ee
where
\[ \Phi_{\De}^2(a,b)=\De(ab)-\De(a)b-(-1)^{|a|}a\De(b). \]
\vsp

We will later specialize to the case of
\vsp

{\bf Definition 4.2} A {\em vertex operator Batalin-Vilkovisky algebra (VOBVA)}
is a GBVA $(\A,1_{\A},\De)$ where $\A$ is a VOSA with the Wick product,
$1_{\A}$ is the vacuum $\va$, and $\De$ is an odd, square zero mode $u_1$.
\vsp

{\bf Remark 4.3} {\em A supercommutative associative algebra is a trivial
example of a
VOBVA, since it is already a VOSA, as explained in} {\bf [LZ2]}.
\vsp

Some of our examples will be differential BV algebras:
\vsp

{\bf Definition 4.4} A {\em differential Batalin-Vilkovisky algebra (DBVA)}
$(\A ,1_{\A},\De ,D,L)$ consists of
\vsp

(i) a superalgebra $\A$ with identity $1_{\A}$,

(ii) an odd, square zero, second order differential operator $\De$ 
($\De 1_{\A}=0$) making $(\A ,\De)$ into a generalized BV algebra,

(iii) an odd, square zero operator $D$ ($D 1_{\A}=0$) and a 
diagonalizable operator $L$ on $\A$ such that
\be {[D,\De ]}=D\De +\De D=L, \label{comm} \ee
and such that the cohomology $\h =$Ker$(D)/$Im$(D)$ of the complex $(\A ,D)$ 
is also a 
GBVA, with the induced product and BV operator $\De$. There are no 
restrictions on the eigenvalues of $L$, which we'll call {\em weights}. In 
most examples $D$ and $\De$ have zero weight as operators, but we'll assume 
only that they are homogeneous (with opposite weights if $L_0\neq 0$). 
\vsp

{\bf Remark 4.5 (i)} {\em One usually thinks of $L$ as the Virasoro operator 
$L_0$, but sometimes we'll just take $L\equiv 0$.}

{\bf (ii)} {\em After finishing the preliminary manuscript I became
aware of Getzler's paper} {\bf [Get2]} {\em in which he defines
a differential BV algebra in a similar fashion. In fact, I was persuaded
by Jim Stasheff that this is a better name than ``cohomology BV algebra''!}
\vsp

{\bf Lemma 4.6} {\em In a VOSA $V$, $L_0-(n+1)\cdot id$ is a derivation of 
the product $\times_n$ defined in Eq.~(\ref{m}). In particular,
$L_0$ is a derivation with respect to the Wick product.}
\vsp

{\it Proof.} \vsp

If $a(z)=\sum a_nz^{-n-1}$ is an $L_0$-homogeneous element with 
$L_0(a)=\mbox{wt}_a\cdot a$, then it follows from the axioms of a VOSA
that
\[ {[L_0,a_n]}=(\mbox{wt}_a-n-1)\cdot a_n \;\;\; \mbox{([FHL])}. \]
Then for homogeneous $a$, $b$,
\bea (L_0-(n+1)\cdot id)(a_n(b))=&&
L_0(a_n(b))-(n+1)a_n(b)\nn \\
=&&(L_0a_n)(b)-(n+1)a_n(b)\nn \\
=&&(a_nL_0)(b)+[L_0,a_n](b)-(n+1)a_n(b)\nn \\
=&&a_n(L_0(b))+(\mbox{wt}_a-n-1)a_n(b)-(n+1)a_n(b)\nn \\
=&&(\mbox{wt}_b+\mbox{wt}_a-n-1-n-1)a_n(b)\nn \\
=&&(\mbox{wt}_a+\mbox{wt}_b-2n-2)a_n(b)\nn \eea
and
\bea &&((L_0-(n+1)\cdot id)(a))_n(b)+a_n((L_0-(n+1)\cdot id)(b))\nn \\
=&&(\mbox{wt}_a-n-1)a_n(b)+(\mbox{wt}_b-n-1)a_n(b)\nn \\
=&&(\mbox{wt}_a+\mbox{wt}_b-2n-2)a_n(b).\nn \eea
Extend to $V$ by bilinearity. $\Box$
\vsp

{\bf Proposition 4.7} {\em In a DBVA we have 
\[ {[L,D]}=[L,\De ]=0, \]
and $\De$ descends to the cohomology $\h$, which has zero $L$-weight.}
\vsp

{\it Proof.}\vsp

${[L,D]}=[\De D+D\De ,D]=\De D^2+D\De D-D\De D-D^2\De =0.$
\vsp

Similarly ${[L,\De ]}=0$. Then both $D$ and $\De$ preserve eigenspaces
of $L$, and anticommute on the kernel. 
Moreover, by virtue of (\ref{comm}) the cohomology can only
be realized in weight zero: if $a$ is a homogeneous element with $L(a)
=\lambda a$ and $D(a)=0$, then
\[ D(\De(a))+\De(D(a))=L(a) \;\;\; \Ra \;\;\; D(\De(a))=\lambda a. \]
If $\lambda\neq 0$, we have $a\in Im(D)$ and $a=0$ in $\h$. $\Box$ 
\vsp

Let us discuss the properties of the generalized BV bracket. Parts of the 
following
computations have appeared in other forms in {\bf [Get1]}, {\bf [Ko]}, 
{\bf [LZ1]}, {\bf [PS]}. We define a new ``Lie'' bracket
on $\A$ which is bilinear and skew-symmetric: Let
\[ {[a,b]}=ab-(-1)^{|a||b|}ba \label{lie} \]
for homogeneous $a$, $b$. Let $\Phi_{\De}^r$ denote the usual $r$-form and
$\tilde{\Phi}_{\De}^r$ the one defined with respect to the bracket above.
This bilinear product was also mentioned in {\bf [LZ1]}. Note that if
$\A$ is supercommutative and associative then ${[\; \; , \; \; ]}$ is trivial,
and so are all $\tilde{\Phi}_{\De}^r$. Finally, recall that $\hat{\A}
=\op_j \A^{j+1}$.
\vsp

{\bf Proposition 4.8} {\it Let $(\A ,\De)$ be a generalized BV algebra. The 
generalized BV bracket satisfies
\vsp

(i) Modified skew-symmetry in $(\hat{\A},\{ \; , \; \})$:
\[ \{ a,b \} +(-1)^{(|a|+1)(|b|+1)}\{ b,a \} =(-1)^{|a|}
\tilde{\Phi}_{\De}^2(a,b). \]

(ii) Leibniz rule (as opposed to the cyclically symmetric version
of the Jacobi identity) in $(\hat{\A},\{ \; , \; \})$:
\[ \{ x, \{ a,b \} \} =\{ \{ x,a \} ,b\} +(-1)^{(|x|+1)(|a|+1)}
\{ a, \{ x,b \} \} . \]

(iii) Poisson rule (superderivation property of $\{ a,\un{{}} \}$ 
on $\A$):
\bea && \{ a,bc \} -\{ a,b \} c-(-1)^{(|a|+1)|b|}b\{ a,c \} \nn \\ &&
=(-1)^{|a|}\Phi_{\De}^3(a,b,c)\equiv 0.\nn \eea

(iv) Derivation rule for $\De$ on $(\hat{\A} ,\{ \; , \; \})$:
\[ \De (\{ a,b \} )-\{ \De(a),b\} -(-1)^{|a|+1}\{ a,\De(b)\}=0. \] }
\vsp

{\bf Proposition 4.9} {\it If $(\A,1_{\A},\De,D,L)$ is a differential BV 
algebra where both $D$ and $L$ are derivations of $\A$, then just like 
$\De$, the operator $D$ is a derivation of $(\hat{\A},\{ \; ,\; \})$. Namely,
\[ D(\{ a,b\}) =\{ D(a),b\} +(-1)^{|a|+1}\{ a,D(b)\}. \]}

{\bf Remark 4.10} {\it This result was proven in {\bf [LZ1]} in the special
case of a BRST complex with $\De =b_0$ ($b(z)=\sum b_nz^{-n-2}$) and
$D=Q=$ BRST operator.}
\vsp

{\bf Remark 4.11} {\it A variant of property~(i) in Proposition~4.8 was written
in {\bf [LZ1]} in the form
\bea &&\{ u,v\} +(-1)^{(|u|-1)(|v|-1)}\{ v,u\} \label{sk} \\
&&=(-1)^{|u|-1}(Qm'(u,v)-m'(Qu,v)-(-1)^{|u|}m'(u,Qv))\nn \eea
(Eq. (2.22)) where $m'$ is a bilinear operation. Since $Q$ is a derivation of
the Wick product, we are justified in calling~(\ref{sk}) ``skew-symmetry up
to homotopy''. But identity~(i) is valid in any GBVA, and we need not have 
the additional DBVA structure. Although 
$\De$ is a derivation of the bracket it defines, we will not go so far
as to say GBVA's are {\em homotopy BV algebras (HBVA)}, because
$\De$ is not a derivation of the original product on $\A$; but
a DBVA is a step towards an HBVA. One could make up several definitions
for an HBVA depending on which identities (including supercommutativity,
associativity, and $\De^2=0$) hold on the nose and which up to homotopy
and higher homotopies. As was noted in
{\bf [LZ1]}, algebraic identities holding only up to homotopy
are discussed in Stasheff's work {\bf [St1]}. Also see 
{\bf [LS]}, {\bf [LZ1]}, {\bf [Get1]}, and {\bf [GV]}.}
\vsp

Lie algebras and generalized BV algebras (with the generalized BV bracket
and shifted grading) are two large classes of examples of what J.-L.~Loday 
calls Leibniz algebras:
\vsp

{\bf Definition 4.12 [Lo]} A {\it Leibniz algebra} $\A =\op_j\A^j$ is a vector
space with super $\Z$-grading and a bilinear superbracket
\[ {[\; , \; ]}:\A \times \A \ra \A \]
satisfying
\[ {[x,[y,z]]}=[[x,y],z]+(-1)^{|x||y|}[y,[x,z]] \]
for all homogeneous $x$, $y$.
\vsp

Then a Leibniz algebra module $M$ is just a Lie algebra module over the
Lie algebra $\bar{\A}$, which is the quotient of $\A$ by the two-sided
ideal generated by all $[x,x]$, $x\in\A$. The classical (Chevalley-Eilenberg)
homology complex can now be replaced by the complex
\[ (T\A \ot M, \del) \]
where $T\A$ is the tensor algebra on $\A$ and $\del$ is an operator which
has almost the same formula as the classical one. Showing $\del$ is square
zero requires nothing more than the Leibniz property of the bracket in $\A$.
This variant of Lie algebra (co)homology is explained in {\bf [Lo]}.
\vsp

For an interesting Leibniz algebra which is not a Lie algebra, see 
{\bf [Mi1]}. Michor defines an odd generalization $\{ \; ,\;\}^2$
of the Poisson bracket to $\Om^{\ast}(M)$ from $\Om^0(M)=C^{\infty}(M)$, where
$M$ is a symplectic manifold. The new bracket satisfies the Leibniz 
identity, but not skew-symmetry (Lemma~4.3). One has to take the quotient
of $\Om^{\ast}(M)$ by the subspace $B^{\ast}(M)$ of exact forms in order to 
restore skew-symmetry.
\vsp

Here are the proofs of the two main propositions.
\vsp

{\em Proof of Proposition 4.8} 
\vsp

(i) $\{ a,b \} +(-1)^{(|a|+1)(|b|+1)}\{ b,a \}$
\vsp

$=(-1)^{|a|}\Phi_{\De}^2(a,b)+(-1)^{|a||b|+|a|+1}\Phi_{\De}^2(b,a)$
\vsp

$=(-1)^{|a|}(\De(ab)-\De(a)b-(-1)^{|a|}a\De(b))
+(-1)^{|a||b|+|a|+1}(\De(ba)-\De(b)a-(-1)^{|b|}b\De(a))$
\vsp

$=(-1)^{|a|}\De(ab-(-1)^{|a||b|}ba)
-(-1)^{|a|}(\De(a)b-(-1)^{|b|(|a|+1)}b\De(a))$

$-(-1)^{|a|}(-1)^{|a|}(a\De(b)-(-1)^{|a|(|b|+1)}\De(b)a)$
\vsp

$=(-1)^{|a|}(\De([a,b])-[\De(a),b]-(-1)^{|a|}[a,\De(b)])$
\vsp

$=(-1)^{|a|}\tilde{\Phi}_{\De}^2(a,b).$
\vsp

(iii) $(-1)^{|a|}(\{ a,bc\} -\{ a,b\} c-(-1)^{|b|(|a|+1)}b \{ a,c\})$
\vsp

$=\Phi_{\De}^2(a,bc)-\Phi_{\De}^2(a,b)c-(-1)^{|b|(|a|+1)}b
\Phi_{\De}^2(a,c)$
\vsp

$=\Phi_{\De}^3(a,b,c)\equiv 0.$
\vsp

(iv) $(-1)^{|a|}(\De(\{ a,b \})-\{ \De(a),b\} -(-1)^{|a|+1}\{ a,\De(b)\})$
\vsp

$=\De(\Phi_{\De}^2(a,b))+\Phi_{\De}^2(\De(a),b)+(-1)^{|a|}\Phi_{\De}^2
(a,\De(b))$
\vsp

$=\De(\De(ab)-\De(a)b-(-1)^{|a|}a\De(b))
+\De(\De(a)b)-\De(\De(a))b-(-1)^{|a|+1}\De(a)\De(b)$

$+(-1)^{|a|}\De(a\De(b))-(-1)^{|a|}\De(a)\De(b)-a\De(\De(b))$
\vsp

$=0.$
\vsp

(ii) We will make use of
\vsp

(A) $\{ x,ab\} -\{ x,a\} b-(-1)^{(|x|+1)|a|}a\{ x,b\} =0.$
\vsp

(B) $\{ x,\De(a)b\} -\{ x,\De(a)\} b-(-1)^{(|x|+1)(|a|+1)}\De(a)\{ 
x,b \}=0.$
\vsp

(C) $\{ x,a\De(b)\}-\{ x,a\} \De(b)-(-1)^{(|x|+1)|a|}a\{ x,\De(b)\}=0.$
\vsp

(D) $\De(\{ x,a\} )-\{ \De(x),a\} -(-1)^{|x|+1}\{ x,\De(a)\} =0.$
\vsp

(E) $\De(\{ x,b\} )-\{ \De(x),b\} -(-1)^{|x|+1}\{ x,\De(b)\} =0.$
\vsp

(F) $\De(\{ x,ab\} )=\{ \De(x),ab\} +(-1)^{|x|+1}\{ x,\De(ab)\}$
\vsp

$=\{ \De(x),a\} b+(-1)^{|a||x|}a\{\De(x),b\}+(-1)^{|x|+1}
\{ x,\De(ab)\} .$
\vsp

Now:
\vsp

$\{ \{ x,a\} ,b\} +(-1)^{(|x|+1)(|a|+1)}\{ a,\{ x,b\} \}$
\vsp

$=(-1)^{|a|+|x|+1}\Phi_{\De}^2(\{ x,a \} ,b) 
+(-1)^{|a||x|+|x|+1}\Phi_{\De}^2(a,\{ x,b\} )$
\vsp

$=(-1)^{|a|+|x|+1}(\De(\{ x,a\} b)-\De(\{ x,a\} )b-(-1)^{|a|+|x|+1}
\{ x,a\} \De(b))$

$+(-1)^{|a||x|+|x|+1}(\De(a\{ x,b\} )-\De(a)\{ x,b\} -(-1)^{|a|}
\De(\{ x,b \} ))$
\vsp

$=(-1)^{|a|+|x|+1}\De(\{ x,a\} b)+(-1)^{|a|+|x|}\De(\{ x,a\}) b
-\{ x,a\} \De(b)$

$+(-1)^{|a||x|+|x|+1}\De(a\{ x,b\} )+(-1)^{|a||x|+|x|}\De(a)\{ x,b\}
+(-1)^{|a||x|+|a|+|x|}a\De(\{ x,b\} )$
\vsp

$=(-1)^{|a|+|x|+1}\De(\{ x,ab\} +(-1)^{|a||x|+|a|+1}a\{ x,b\} )$ by (A)

$+(-1)^{|a|+|x|}(\{ \De(x),a\} +(-1)^{|x|+1}\{ x,\De(a)\} )b$ by (D)

$-\{ x,a\} \De(b)+(-1)^{|a||x|+|x|+1}\De(a\{ x,b\} ) 
+(-1)^{|a||x|+|x|}\De(a)\{ x,b\}$

$+(-1)^{|a||x|+|a|+|x|}a(\{ \De(x),b\} +(-1)^{|x|+1}\{ x,\De(b)\} )$ by (E)
\vsp

$=(-1)^{|a|+|x|+1}\De(\{ x,ab\} )+(-1)^{|a|+|x|}\{ \De(x),a\} b
+(-1)^{|a|+1}\{ x,\De(a)\} b-\{ x,a\} \De(b)$

$+(-1)^{|a||x|+|x|}\De(a)\{ x,b\}
+(-1)^{|a||x|+|a|+|x|}a\{ \De(x),b\}
+(-1)^{|a||x|+|a|+1}a\{ x,\De(b)\}$
\vsp

$=(-1)^{|a|+|x|+1}\De(\{ x,ab\})
+(-1)^{|a|+|x|}\{ \De(x),a\} b$

$+((-1)^{|a|+1}\{ x,\De(a)\} b+(-1)^{|a||x|+|x|}\De(a)\{ x,b\} )$

$+(-\{ x,a\} \De(b)+(-1)^{|a||x|+|a|+1}a\{ x,\De(b)\} )
+(-1)^{|a||x|+|a|+|x|}a\{ \De(x),b\}$
\vsp

$=(-1)^{|a|+|x|+1}\De(\{ x,ab\} )+(-1)^{|a|+|x|}\{ \De(x),a\} b
+(-1)^{|a||x|+|a|+|x|}a\{ \De(x),b\}$

$+(-1)^{|a|+1}\{ x,\De(a)b\}$ by (B)

$-\{ x,a\De(b)\}$ by (C)
\vsp

$=(-1)^{|a|+|x|}(-\De(\{ x,ab\} )+\{ \De(x),a\} b+(-1)^{|a||x|}a
\{ \De(x),b\} )$

$+(-1)^{|a|+1}\{ x,\De(a)b\} -\{ x,a\De(b)\}$
\vsp

$=(-1)^{|a|}\{ x,\De(ab)\}$ by (F)

$+(-1)^{|a|+1}\{ x,\De(a)b\} -\{ x,a\De(b)\}$
\vsp

$=(-1)^{|a|}\{ x,\De(ab)-\De(a)b-(-1)^{|a|}a\De(b)\}$
\vsp

$=\{ x,(-1)^{|a|}\Phi_{\De}^2(a,b)\}$
\vsp

$=\{ x,\{ a,b\} \} . \; \; \; \Box$
\vsp

{\it Proof of Proposition 4.9.}
\vsp

$D(\{ a,b\})-\{ D(a),b\} -(-1)^{|a|+1}\{ a,D(b)\}$
\vsp

$=(-1)^{|a|}D(\De(ab)-\De(a)b-(-1)^{|a|}a\De(b))$

$+(-1)^{|a|}(\De(D(a)b)-(\De D)(a)b+(-1)^{|a|}D(a)\De(b))$

$+\De(aD(b))-\De(a)D(b)-(-1)^{|a|}a(\De D)(b)$
\vsp

$=(-1)^{|a|}(D\De)(ab)-(-1)^{|a|}(D\De)(a)b+\De(a)D(b)$

$-D(a)\De(b)-(-1)^{|a|}a(D\De)(b)+(-1)^{|a|}\De(D(a)b)$

$-(-1)^{|a|}(\De D)(a)b+D(a)\De(b)+\De(aD(b))$

$-\De(a)D(b)-(-1)^{|a|}a(\De D)(b)$
\vsp

(3rd, 4th, 8th, 10th terms cancel out. Expand 1st term, join 2nd and
7th, also 5th and 11th, and 6th and 9th.)
\vsp

$=-(-1)^{|a|}\De(D(ab))+(-1)^{|a|}L(ab)-(-1)^{|a|}(D\De +\De D)(a)b$

$-(-1)^{|a|}a(D\De +\De D)(b)+(-1)^{|a|}\De(D(a)b+(-1)^{|a|}aD(b))$
\vsp

(Join 1st and 5th, also 2nd, 3rd, and 4th terms.)
\vsp

$=-(-1)^{|a|}\De(D(ab)-D(a)b-(-1)^{|a|}aD(b))
+(-1)^{|a|}(L(ab)-L(a)b-aL(b))$
\vsp

$=-(-1)^{|a|}\De(\Phi_D^2(a,b))+(-1)^{|a|}(\Phi_L^2(a,b))=0.\;\;\;\Box$
\vsp

Although the level of generalization in this section seems adequate for
most examples, it is instructive to study the bracket $\{ a,b\}$ on a
superalgebra $\A$ with a plain odd operator $\De$, not necessarily 
square zero, or of order two. We find that 
Proposition~4.8 is replaced by
\vsp

{\bf Proposition 4.13} {\it For a superalgebra $\A$ and an odd operator
$\De$, the bracket $\{ \;\; ,\;\; \}$ defined by (\ref{bracket}) satisfies
the following identities:}

{\it (i) Modified skew-symmetry:}
\[ \{ a,b\} +(-1)^{(|a|+1)(|b|+1)}\{ b,a\} =(-1)^{|a|}\tilde{\Phi}_{\De}^2
(a,b).\]

{\it (ii) Modified Leibniz rule:}
\bea \{\{ x,a\} ,b\} +&&(-1)^{(|x|+1)(|a|+1)}\{ a,\{ x,b\}\} -\{ x,
\{ a,b\}\}\nn \\ =&&(-1)^{|a|}( \De(\Phi_{\De}^3(x,a,b))-\Phi_{\De^2}^3
(x,a,b)+\Phi_{\De}^3(\De(x),a,b)\nn\\ &&+(-1)^{|x|}\Phi_{\De}^3(x,\De(a),
b)+(-1)^{|x|+|a|}\Phi_{\De}^3(x,a,\De(b)) ).\nn \eea

{\it (iii) Modified Poisson rule:}
\[ \{ a,bc\} -\{ a,b\} c-(-1)^{(|a|+1)|b|}b\{ a,c\} =(-1)^{|a|}
\Phi_{\De}^3(a,b,c).\]

{\it (iv) Modified derivation rule for $\De$:}
\[ \De(\{ a,b\} )-\{ \De(a),b\} -(-1)^{|a|+1}\{ a,\De(b)\}=(-1)^{|a|}
\Phi_{\De^2}^2(a,b).\]

It is clear that we obtain a classical BV algebra when $\A$ is 
supercommutative, $\De$ is of order two, and $\De^2$ is zero (or
of order one). The associativity condition on classical BV algebras 
turns out to be superfluous. The proof of Proposition~4.13 can be obtained
from that of Proposition~4.8 mostly by retaining terms (such as those
containing $\De^2$ and $\Phi_{\De}^3$) which were formerly discarded.
Part~(ii) has a strong resemblance to Lemma~1.5 in {\bf [Ko]}, and
(iv) is exactly his Equation~(1.8).
\vsp

If we now add an odd, square zero, first order differential operator
$D$ and a diagonalizable, first order differential operator $L$ to
this picture, such that $D\De +\De D=L$ as before, Proposition~4.9
remains unchanged. The result is already in terms of $\De(\Phi_{D}^2
(a,b))$ and $\Phi_{L}^2(a,b)$, allowing more general operators. The fact
that $D^2=0$ is not used at all. The need for unadorned algebras and
operators will become clear in Section~4.2.2, where we introduce
topological chiral algebras.
\vsp

\vsp

\vsp

{\it 4.2 Examples of Vertex Operator BV Algebras}
\vsp

{\it 4.2.1 The Vertex Operator Weil Algebra}
\vsp

The semi-infinite Weil complex $W^{\infty /2}\g$ associated to a {\em tame}
Lie algebra $\g$ (i.e. $\g =\oplus_n\g_n$, $[\g_m,\g_n]\subset\g_{m+n}$, 
dim$\g_n<\infty$) was first considered by Feigin and Frenkel in {\bf [FF]}.
They in particular computed the semi-infinite cohomology for $\g =Witt$.
Next, the {\em vertex operator Weil algebra (VOWA)} 
\[ \W =S^{\infty /2}\tilde{\ell}'\ot_{\Co}
\wedge^{\infty /2}\tilde{\ell}' \]
on a {\em loop algebra}
\[ \tilde{\ell}=\ell\ot_{\Co}\Co [t,t^{-1}] \]
with bracket
\[ {[x\ot t^m,y\ot t^n]}=[x,y]\ot t^{m+n} \]
(where $\ell$ is a finite dimensional complex Lie algebra) was studied
by the author in {\bf [A2, A3]}. We refer the reader to {\bf [A3]}
for details.
\vsp

The VOWA is graded by the eigenvalues 0, 1, 2, ... of $L_0$, and 
$W^{\infty /2}\tilde{\ell}[0]$ is exactly the classical Weil algebra
\[ W\ell =S\ell'\ot\wedge\ell' \]
considered in Section 3.2.3. Most semi-infinite operators restrict to the
classical subspace. In {\bf [A3]} it was noted that $\W$ is a 
topological chiral algebra (see Section 4.2.2) because of an identity
\[ hk+kh=-L_0 \]
involving the semi-infinite Koszul differential, a new semi-infinite
homotopy operator, and the Virasoro operator $L_0$.
If we don't mind taking $L=0$ in (\ref{comm}), we have several VOBVA and 
DBVA structures on $\A =\W$ ($Q$ is the semi-infinite cohomology operator):
\vsp

(i) $(\A ,\De ,D,L)=(\W ,k,h,-L_0),\;\;\; (\h ,\De)=(\Co ,0)$
\vsp

(ii) $(\A ,\De ,D,L)=(\W ,Q,h,0),\;\;\; (\h ,\De)=(\Co ,0)$
\vsp

(iii) $(\A ,\De ,D,L)=(\W ,h,Q,0),\;\;\; (\h ,\De)=(H(\W ,Q),h)$
\vsp

(iv) $(\A ,\De ,D,L)=(\W ,k,Q,0),\;\;\; (\h ,\De)=(H(\W ,Q),k).$
\vsp

Now, the space $\h =H(\W ,Q)$ is not a supercommutative, associative
algebra, but it contains the classical space $\h [0]=H(W\ell ,d)$, which 
{\sl is}. As was shown in {\bf [A2]}, $\h$ is infinitely
richer than $H(W\ell ,d)$. The operators $h$ and $k$ induce nontrivial 
first and second order differential operators on $\h$, and case~(iv) even has 
a (possibly) nontrivial BV bracket on the cohomology.
For reductive $\ell$, two more vertex operators
and their weight zero modes $r$, $t$ can be exploited {\bf [A2, A3]}.
These additional operators also satisfy
\[ rt+tr=-L_0. \]
We can now add the following to our list of differential BV algebras:
\vsp

(v) $(\A ,\De ,D,L)=(\W ,t,r,-L_0),\;\;\; (\h ,\De)=(H(W\ell ,r),t)$
\vsp

(vi) $(\A ,\De ,D,L)=(\W ,Q,r,0),\;\;\; (\h ,\De)=(H(W\ell ,r),Q)$
\vsp

(vii) $(\A ,\De ,D,L)=(\W ,r,Q,0),\;\;\; (\h ,\De)=(H(\W ,Q),r)$
\vsp

(viii) $(\A ,\De ,D,L)=(\W ,t,Q,0),\;\;\; (\h ,\De)=(H(\W ,Q),t)$
\vsp

(once again, (v) and (viii) have a potentially nonzero BV bracket as $t$ is
a genuine second order operator) and also
\vsp

(ix) $(\A ,\De ,D,L)=(\W ,r,h,0),\;\;\; (\h ,\De)=(\Co ,0)$
\vsp

(x) $(\A ,\De ,D,L)=(\W ,t,h,0),\;\;\; (\h ,\De)=(\Co ,0)$
\vsp

(xi) $(\A ,\De ,D,L)=(\W ,h,r,0),\;\;\; (\h ,\De)=(H(W\ell ,r),h)$
\vsp

(xii) $(\A ,\De ,D,L)=(\W ,k,r,0),\;\;\; (\h ,\De)=(H(W\ell ,r),k).$
\vsp

While $h$ has a classical analogue,
$k$, $r$, and $t$ don't, and as noted above, $k$ and $t$ may give
rise to nonzero BV brackets in the cohomology. So far
we have been using a restricted version of Definition~4.1 where the operator 
$D$ is a derivation, which is very common. Relaxing this condition (say, 
allowing $D$ to be a second order differential operator),
we would get several additional examples. One of the DBVA's
we can obtain from $\W$ in this fashion is
\vsp

(xiii) $(\A,\De,D,L)=(\W,Q,k,0),\;\;\; (\h,\De)=(W\ell ,d)$,
\vsp

so that the classical Weil algebra of 3.2.3 is the cohomology of a BV algebra
(with trivial BV bracket).
\vsp

\vsp

\vsp

{\it 4.2.2 Topological Chiral Algebras and the BRST Complex}
\vsp

We take the following definition from {\bf [LZ1]}. A {\em topological
chiral algebra (TCA)} consists of
\vsp

(i) A VOSA,

(ii) A weight one even field $F(z)=\sum F_nz^{-n-1}$ whose residue
(charge) $F_0$ is the ``fermion number operator'', or ``ghost number
operator'',

(iii) A weight one primary (Virasoro-singular) field $J(z)=\sum J_n
z^{-n-1}$ with fermion number one and a square zero charge $Q=J_0$,

(iv) A weight two primary field $G(z)=\sum G_nz^{-n-2}$ with fermion
number $-1$, satisfying
\be {[Q,G(z)]}=L(z),\label{c1}\ee
where $L(z)=\sum L_nz^{-n-2}$ is the ``stress-energy (Virasoro) field''.
By (\ref{c1}) the Virasoro algebra acts trivially on the cohomology, hence
the name ``topological''.
\vsp

The above definition is tailored for the {\em BRST complex} 
\be V\ot\wedge^{\ast} \label{c2}\ee
where V is any VOA with central charge 26 and $\wedge^{\ast}$ is the
``ghost system'', i.e. a simple $bc$-system generated by two fields
\[ b(z)=\sum b_nz^{-n-2},\;\;\; c(z)=\sum c_nz^{-n+1}\]
with stress-energy field
\[ L^{\wedge}(z)=-:(\frac{d}{dz}b(z))c(z):-2:b(z)\frac{d}{dz}c(z):\; .\]
{}From the viewpoint of semi-infinite cohomology, (\ref{c2}) is a cohomology
complex where
\be \wedge^{\ast}=\wedge^{\infty /2+\ast}(Witt)',\label{c3}\ee
the Witt algebra being the quotient of $Vir$ by its center. The central
element of $Vir$ acts by $-26$ on~(\ref{c3}), so in order to have a closed 
action of $Witt$ and a square zero differential we have to tensor with
a VOA of central charge 26 (``cancellation of anomalies''). The choice of
vacuum is given by
\[ b_n\va =c_m\va =0\;\;\; \mbox{for $n\geq -1,m\geq 2$,}\]
because (\ref{c3}) is a VOSA only in this vacuum. The Virasoro (Witt)
action is given by
\[ \rho(z)=L^{\wedge}(z) \;\;\;\mbox{on $\wedge^{\ast}$}\]
and
\[ \theta(z)=L(z)=L^{V}(z)+L^{\wedge}(z)\;\;\;\mbox{on $V\ot\wedge^{\ast}$.}
\]
We also have 
\bea &&\io(L_n)=b_n=G_n, \;\;\;\ep(L_n')=c_{-n},\;\;\;
d=Q=J_0,\nn \\ &&J(z)=:c(z)(L^{V}(z)+\frac{1}{2}L^{\wedge}(z)):,\nn \eea
so that (\ref{c1}) is just the Cartan identity in disguise. The same equation
tells us that the vertex operator $L(z)$ is $Q$-exact (not true, for example,
for the VOWA). Finally,
\[ F(z)=:c(z)b(z):\; .\]
The vertex operator Weil algebra is also a TCA and a DBVA with 
\[ {[h,k(z)]}=-L(z)\]
and
\[ F(z)=\sum_u :c^{u'}(z)b^u(z): \]
(see {\bf [A3]}). If $\ell$ is semisimple, it is also true that
\[ {[r,t(z)]}=-L(z).\]
\vsp

Topological chiral algebras where $G_0^2=0$ are both VOBVA's and DBVA's, with
\be \De =Q,\;\;\; D=G_0,\;\;\; L=L_0\label{c4}\ee
or
\[ \De =G_0,\;\;\; D=Q, \;\;\; L=L_0.\]
In the case of a BRST complex, (\ref{c4}) yields zero cohomology for
$D=b_0$ (because $b_0c_0+c_0b_0=1$), and in any case the BV bracket on the
complex is trivial. The second choice is the DBVA which was implicit in
{\bf [LZ1]}.
\vsp

There seems to be a demand for TCA's without the condition $G_0^2=0$
as well, as was explained to the author by Jos\'{e}~M.~Figueroa-O'Farrill
after the first preprint of this article appeared (see {\bf [Ka, Get2, IR,
Fi, FS]}). It is once again desirable to have a classical BV algebra
as the $Q$-cohomology. Even without the ``primary'' condition on $G(z)$
(see (\ref{c1})),
one can show that Proposition~4.7 still holds for a TCA: We have
\[ {[L_0,Q]}=[L_0,G_0]=0\]
because $Q$ and $G_0$ are weight zero modes, and all the arguments in the
proof are valid. Then the cohomology ${\cal H}$, being a VOSA of weight
zero, is a supercommutative, associative algebra with a second order,
odd differential operator $G_0$. The only remaining issue is to impose 
minimal conditions on $G(z)$ so that $G_0^2=0$ in ${\cal H}$. In {\bf [Ka]}
and {\bf [Get2]} we find {\it Kazama algebras} in which $G_0$ satisfies
an identity
\[ G_0^2=[Q,P]\]
where $P$ is an odd, third order operator, making $G_0^2=0$ in ${\cal H}$.
We observe that the above condition will be satisfied if $G(z)$ is primary:
\vsp

{\bf Proposition 4.14} {\it Let $V$ be a topological chiral algebra with
a primary field $G(z)$. Then ${\cal H}=H^{\ast}(V,Q)$ is a classical BV 
algebra.}
\vsp

{\it Proof.} Since $G(z)$ is a primary field of weight two, we have the
identity
\be {[L_m,G_n]}=(m-n)G_{m+n}\label{primary}\ee
(see {\bf [FLM]}), in particular
\[ {[L_{-2},G_{-2}]}=0.\]
Then the field $:G(z)G(z):$, or the associated element $(G_{-2})^2\vac$,
is a cocycle:
\[ Q\; (G_{-2})^2\vac =[Q,G_{-2}]G_{-2}\vac -G_{-2}[Q,G_{-2}]\vac
=L_{-2}G_{-2}\vac -G_{-2}L_{-2}\vac =0.\]
Also, since $Q$ commutes with all $L_n$, we have
\[ Q\; L_3(G_{-2})^2\vac =0.\]
But then by the identity
\[ QG_0+G_0Q=L_0,\]
both $(G_{-2})^2\vac$ and $L_{3}(G_{-2})^2\vac$ (of weights 4 and 1
respectively) are coboundaries. Recall that since $Q$ is a residue,
every coboundary vertex operator $Q\; R(z)$ is a $Q$-commutator,
namely 
\[ Q\; R(z)=[Q,R(z)]=\sum [Q,R_{(n)}]z^{-n-1}.\]
The proof will be complete when we show that $G_0^2$ is a linear 
combination of the weight zero modes of the fields corresponding to
$(G_{-2})^2\vac$ and $L_3(G_{-2})^2\vac$. By identity~(\ref{geb1})
we find that the weight zero mode of $:G(z)G(z):$ is
\bea (G_{-2}G)_0=(G_{(-1)})_{(3)} &&=\sum_{i\geq 0}(-1)^i(-1)^i
(G_{(-1-i)}G_{(3+i)}-G_{(2-i)}G_{(i)})\nn\\
&&=-G_{(1)}G_{(1)}-G_{(2)}G_{(0)}-G_{(0)}G_{(2)}=-G_0^2-[G_1,G_{-1}],\nn\eea
where we translate between the two conventions in subscripts by
\[ R(z)=\sum R_nz^{-n-wt(R)}=\sum R_{(n)}z^{-n-1}.\]
Next, we observe that 
\[ L_3(G_{-2})^2\vac =[L_3,G_{-2}]G_{-2}\vac +G_{-2}[L_3,G_{-2}]\vac
=5G_1G_{-2}\vac +5G_{-2}G_1\vac =5G_1G_{-2}\vac ,\]
by (\ref{primary}) and the fact that coefficients of $z^n$ for $n<0$
annihilate the vacuum. We calculate the zero mode of $\frac{1}{5}L_3
(G_{-2})^2\vac$ as follows:
\[ (G_1G)_0=(G_{(2)}G)_{(0)}=2(-2G_{(1)}^2+[G_{(2)},G_{(0)}])
=2(-2G_0^2+[G_1,G_{-1}]),\]
and $G_0^2$ is indeed a linear combination of $(G_{-2}G)_0$ and $(G_1G)_0$.
$\Box$
\vsp

This result and several similar ones are reported to be known to 
J.~M. Figueroa-O'Farrill
and T. Kimura. Their proof seems to involve an identity similar 
to~(\ref{geb2}), and a modified version would read as follows:
\[ 2G_0^2=[G_0,G_0]=[G_{(1)},G_{(1)}]=(G_{-1}G)_0+(G_0G)_0\]
by (\ref{geb2}), where once again
\[ Q\; G_{-1}G=Q\; G_0G=0\]
is a consequence of (\ref{primary}).
\vsp

\vsp

\vsp

{\bf 5 BV Master Equation}
\vsp

{\it 5.1 Quantum BV Master Equation for Classical BV Algebras}
\vsp

In {\bf [BV]} it is proposed that we look for an action (a bosonic function
of fields and antifields) 
\be W=S+\sum_{p=1}^{\infty}\hbar^pM_p\label{d1}\ee
satisfying the {\em quantum master equation}
\be \{ W,W\} =2i\hbar\De(W),\label{d2}\ee
or equivalently,
\bea &&\{ S,S\} =0\nn\\
&&\{ M_1,S\} =i\De(S)\label{d3}\\
&&\{ M_p,S\}=i\De(M_{p-1})-\frac{1}{2}\sum_{q=1}^{p-1}\{ M_q,M_{p-q}\}
\;\;\;\mbox{for $p\geq 2$,}\nn\eea
so that
\be \De(exp(\frac{i}{\hbar}W))=0.\label{d4}\ee
Note that the classical part of $W$, shown by $S$, satisfies the {\em 
classical master equation}
\be \{ S,S\} =0.\label{d5}\ee
Infinite sums will be assumed to have a purely formal 
meaning, or else to be finite when interpreted in some context. Our 
task will be to understand why condition~(\ref{d4}) follows from~(\ref{d2})
(a result often cited in literature without proof) when $\A$ is a 
supercommutative, associative algebra.
\vsp

Assume 
\be \{ W,W\} =\lambda\De(W)\label{d6}\ee
for some even element $W$ of a classical BV algebra and some complex constant
$\lambda$. Note that for even $V$, $W$ we have
\[ \{ V,W\} =\Phi_{\De}^2(V,W).\]
\vsp

{\bf Lemma 5.1} {\it For $k\geq 0$,} 
\be \{ W,W^k\} =k\lambda \De(W)W^{k-1}. \label{bik} \ee

{\em Proof.} Induction on $k$. For $k=1$, we have $\{ W,W\} =\lambda\De(W)$.
Assume (\ref{bik}) holds for 1, 2, ..., $k$.
\bea \{ W,W^{k+1}\} &&=\{ W,W\cdot W^k\} =\{ W,W\} W^k+W\{ W,W^k\} \nn \\
&&=\lambda\De(W)W^k+W(k\lambda\De(W)W^{k-1}) \; \; \; 
\mbox{by (\ref{d6}) and (\ref{bik})}\nn \\ &&=(k+1)\lambda\De(W)W^k .
\, \, \, \nn\Box\eea

{\bf Lemma 5.2} {\it For $k\geq 0$,}
\be \De(W^k)= \frac{k(k-1)}{2}\lambda\De(W)W^{k-2}+k\De(W)W^{k-1}.
\label{cik} \ee

{\em Proof.} Induction on $k$. The statement holds for $k=0$, $1$.
Assume (\ref{cik}) holds for 1, 2, ..., $k$.
\bea \De(W^{k+1}) &&= \De(W\cdot W^k)=\Phi_{\De}^2(W,W^k)+\De(W)W^k+
W\De(W^k) \nn \\ &&=\{ W,W^k\} +\De(W)W^k+W\De(W^k)\nn \\ 
&&=k\lambda\De(W)W^{k-1}+\De(W)W^k+W(\frac{k(k-1)}{2}\lambda\De(W)
W^{k-2}+k\De(W)W^{k-1}) \nn \\ &&\mbox{(by Lemma 5.1 and induction step)} \nn
\\ &&=(\frac{k(k-1)}{2}+k)\lambda\De(W)W^{k-1}+(k+1)\De(W)W^k \nn \\
&&=\frac{k(k+1)}{2}\lambda\De(W)W^{k-1}+(k+1)\De(W)W^k. \; \; \; \Box \nn \eea

{\bf Proposition 5.3} $\De(exp(\frac{i}{\hbar}W))=(\frac{i\lambda}{2\hbar}+1)
\De(\frac{i}{\hbar}W)exp(\frac{i}{\hbar}W).$
\vsp

{\em Proof.} We will replace $\frac{i}{\hbar}W$ by $V$, and 
condition~(\ref{d6}) by 
\[ \{ V,V\} =\mu\De(V),\]
where $\mu =i\lambda /\hbar$, and prove the simpler identity
\[ \De(exp(V))=(\frac{\mu}{2}+1)\De(V)exp(V).\]
Thus
\bea &&\De(exp(V)) \nn \\
&&=\De(V)+\sum_{k=2}^{\infty}\frac{k(k-1)}{2k!}\mu\De(V)V^{k-2}
+\sum_{k=2}^{\infty}\frac{k}{k!}\De(V)V^{k-1} \nn \\
&&=\De(V)+\frac{\mu}{2}\De(V)\sum_{k=0}^{\infty}\frac{V^k}{k!}
+\De(V)\sum_{k=1}^{\infty}\frac{V^k}{k!} \nn \\
&&=\frac{\mu}{2}\De(V)\sum_{k=0}^{\infty}\frac{V^k}{k!}+\De(V)
\sum_{k=0}^{\infty}\frac{V^k}{k!}=(\frac{\mu}{2}+1)\De(V)exp(V).
\; \; \; \Box \nn \eea

As a result, we have
\vsp

{\bf Corollary 5.4} {\it If (\ref{d6}) holds for even $W$ with $\lambda 
=2i\hbar$, then 
\[ \De(exp(\frac{i}{\hbar}W))=0.\]}

{\bf Remark 5.5} {\it A recent formulation of Eq.~(\ref{d4}) in terms of the
forms $\Phi_{\De}^k$ was given by J.~Alfaro and P.~H. Damgaard {\bf [AD]}.
This is an identity I completely missed in my preprints, and I am 
grateful to P.H. Damgaard for explaining the equivalence to me. One can show,
from the remarks at the end of Section~2.1 and by induction, that 
\[ \De(W^k)=\sum_{j=1}^k\kj W^{k-j}\Phi_{\De}^j(W,\cdots ,W) \]
for every $k\geq 1$, with the assumptions that $W$ is an even element in 
a supercommutative, associative algebra, and $\De$ is any odd operator. Then
it is a matter of simple algebra to show
\[ \De(exp(\frac{i}{\hbar}W))=exp(\frac{i}{\hbar}W)\sum_{k=1}^{\infty}
\frac{1}{k!}(\frac{i}{\hbar})^k\Phi_{\De}^k(W,\cdots ,W), \]
under the same assumptions. This last identity allows us to 
replace (\ref{d4}) by
\be \sum_{k=1}^{\infty}\frac{1}{k!}(\frac{i}{\hbar})^k\Phi_{\De}^k
(W,\cdots ,W)=0, \ee
which would work very well even for more general types of algebras! Note
that if $\De$ is of order two, then only the first two terms survive,
and we obtain the quantum BV master equation~(\ref{d2}).}
\vsp

\vsp

{\it 5.2 Quantum BV Master Equation for Vertex Operator BV Algebras}
\vsp

A question raised in {\bf [LZ1]} is the meaning of the 
equations~(\ref{d2}) and~(\ref{d5}) in a conformal field theory,
and their relations to deformations of the theory. It is stated that
\vsp

{\bf Proposition 5.6} (Proposition 3.3 in {\bf [LZ1]}) {\em The first
order pole of the operator product expansion (OPE) $(b_{-1}S)(z)(b_{-1}S)(w)$,
where $S$ is an even, weight zero, BRST-invariant element, vanishes
if and only if the bracket $\{ S,S\}$ defined via the BV operator
\[ \De =b_0, \;\;\; \mbox{where}\;\;\; b(z)=\sum b_nz^{-n-2}\]
vanishes.}
\vsp

Such an element will correspond to a ``first order deformation''
(perturbation) of the BRST complex, and if all poles vanish,
one obtains deformations of all orders. The authors then determine
all solutions of $\{ S,S\} =0$ in the BRST cohomology of the ``c=1
model''. 
\vsp

Let us give a proof of Proposition~5.6.
\vsp

{\it Proof.} The n-th order pole of the OPE $u(z)v(z)$ is given by
$u_{n-1}v$ (see {\bf [LZ3]}). Replace $b_n$ by $u_{n+1}$ to achieve
standard notation. We want to prove that if $S$ is even, with 
\be L_0S=QS=0,\label{d15}\ee
then
\[ (u_0S)_0u_0S=0\;\;\;\Lr\;\;\; (u_0S)_0S=0.\]
Equivalently, we want
\[ \{ S,u_0S\} =0\;\;\; \Lr\;\;\; \{ S,S\} =0.\]
One direction is quite generally true:
\[ (u_0S)_0S=0\;\;\;\Ra\;\;\; u_0(u_0S)_0S=0\;\;\;
\Ra\;\;\; -(u_0S)_0u_0S+[u_0,(u_0S)_0]S=0.\]
But
\[ {[u_0,(u_0S)_0]}=[u_0,[u_0,S_0]]=[[u_0,u_0],S_0]-[u_0,[u_0,S_0]]\]
(where $u_0$ is odd, square zero)
\[ \Ra\;\;\; [u_0,[u_0,S_0]]=0\;\;\;\Ra\;\;\; (u_0S)_0u_0S=0.\]
Conversely, assume that
\[ (u_0S)_0u_0S=0,\]
or, as we showed above,
\[ u_0(u_0S)_0S=0.\]
This last equation translates to
\[ b_{-1}\{ S,S\} =0\;\;\;\Ra\;\;\; Qb_{-1}\{ S,S\} =0\;\;\;\Ra\;\;\;
-b_{-1}Q\{ S,S\} +L_{-1}\{ S,S\} =0. \]
But recall that we are working with a CBVA and $Q$ is a derivation.
By virtue of~(\ref{d15}) and Proposition~4.9, we have 
\[ Q\{ S,S\} =0,\]
leading to
\be L_{-1}\{ S,S\} =0.\label{d20}\ee
By a fundamental axiom of VOSA's, condition~(\ref{d20}) means that the
vertex operator $\{ S,S\} (z)$ is a constant, and 
\[ {[ \{ S,S \}_{-1},v(z)]}=0\]
for any $v(z)$, as
\[ \{ S,S\}_iv=0\cdot v=0\;\;\;\forall i\geq 0\]
(see (\ref{geb2})). This in turn implies that 
\[ v_0\{ S,S\} =[v_0,\{ S,S\}_{-1}]\va =0\;\;\;\forall v.\]
In particular for $v_0=F_0$, the ghost number (superdegree) operator,
we have
\[ F_0\{ S,S\} =-\{ S,S\} =0.\;\;\; \Box\]
\vsp

{\bf Remark 5.7} {\it Note that we did not need the condition $L_0S=0$.}
\vsp

Since the proof is valid in any TCA, we may choose one 
with overall nonnegative grading, so that all poles except the first two 
vanish when S is of weight zero.
The VOWA provides examples of such TCA's. The two basic identities are
\[ {[h,k(z)]}=-L(z) \;\;\;\mbox{and}\;\;\; [r,t(z)]=-L(z),\]
the latter being defined for semisimple $\ell$ only.
\vsp

We also make the following observation by comparing weights on both sides
of the master equation:
\vsp

{\bf Proposition 5.8} {\em In a VOBVA, any even solution $W$ of the quantum
master equation~(\ref{d6}) with $\lambda\neq 0$ and $\De(W)\neq 0$ which is
$L_0$-homogeneous has to be of weight zero. In fact, if W is any solution with
a homogeneous component $W'$ of nonzero highest (or lowest) weight, then $W'$
satisfies the classical master equation $\{ W',W'\} =0$.}
\vsp

\vsp

\vsp

{\it 5.3 Quantum BV Master Equation for Generalized BV Algebras}
\vsp

A classical problem: If $(\Li^{\ast},\delta)$ is a $\Z$-graded super Lie
algebra with an odd derivation $\delta$ of the Lie bracket (say of degree~1
and square zero), for which odd elements $a$ of $\Li$ will the addition of
the inner derivation $\ad(a)$ result in an odd, square zero derivation~$\td$?
We first note that $\delta +\ad(a)$ is still a derivation of the Lie bracket.
Next, we observe 
\bea && \td^2=0 \;\;\; \Leftrightarrow \;\;\; \delta^2+[\delta,\ad(a)]
+\ad(a)^2=0 \label{h1} \\ && \Leftrightarrow\;\;\; \ad(\delta(a))+\f{1}{2}\ad
([a,a])=0\;\;\; \Leftrightarrow\;\;\;\ad(\delta(a)+\f{1}{2}[a,a])=0. \nn \eea
This last condition holds when, for example,
\be \delta(a)+\f{1}{2}[a,a]=0. \label{h2}\ee
If ${[\; ,\; ]}$ is an odd Lie bracket, as in a classical BV algebra, we
seek an even element $a$ instead. We may call~(\ref{h2}) a {\it
deformation equation} (see {\bf [NR, Ger]}). In a (classical) differential BV
algebra $(\A ,1_{\A},\De,D,L)$ where $D$ and $L$ are derivations of the
associative product in $\A$, we have {\sl two} deformation equations,
namely
\be D(a)+\f{1}{2}\{ a,a\} =0 \label{h3}\ee
and
\be \De(a)+\f{1}{2}\{ a,a\} =0, \label{h4}\ee
where the second one also deserves the name ``quantum BV master equation''.
Note that although $\De$ is in general a second order differential operator
on $\A$, it is a derivation of $(\hat{\A},\{ \; ,\; \} )$. Then solving the
quantum BV master equation means finding deformations of the BV
operator $\De$ by inner derivations (see  {\bf [St2]}).
Addition of an inner derivation, of course, does not change the BV bracket.
\vsp

In the case of a generalized BV algebra $\A$ with bracket $\{ \; ,\;\}$,
$(\hat{\A},\{ \; ,\;\})$ is only a Leibniz algebra. But all we need 
in~(\ref{h1}) is the Leibniz property and the fact that $\De$ is a derivation
of the bracket, so that the even solutions of Eq.~(\ref{h4}) still give us
deformations of $\De$. In other words, $\De +\{ a,\; \; \}$ is again an odd,
square zero, second order differential operator on $\A$ which induces the
bracket $\{ \;\; ,\;\; \}$ associated with $\De$. If (\ref{h4}) holds in a
DBVA, the commutator of $\De +\{ a,\;\; \}$ with $D$ is still $L$, provided
that $a$ is in the kernel of $D$.
\vsp

\vsp

\vsp

{\it Acknowledgments.} I would like to thank Jean-Louis Koszul for his kind 
remarks and for telling me about Grothendieck's recursive definition.
I am, as always, indebted to Gregg Zuckerman for
his thorough reading of the preliminary versions and for his innumerable
historical remarks, suggestions, and corrections. This article was more or
less deformed into its final shape during a fruitful visit to Yale University,
where Bong Lian also made several contributions involving vertex operator
algebras. I am in particular grateful to G.Z. for bringing the deformation
equation and Loday's book to my attention, and to B.L. for explaining the 
role of the ghost number operator in the proof of Proposition~5.6. After
circulating a preprint, I received considerable encouragement and very
helpful comments from Jim Stasheff, which resulted in the elimination
of many small errors and a marked improvement of the exposition. 
During my visit to the University of North Carolina, he pointed out
connections with several topics (such as strongly homotopy Lie algebras)
and gave me references. Thanks are also due Poul Damgaard for a lengthy
discussion of the quantum BV master equation. Finally, I would like to
thank Jos\'{e} Figueroa-O'Farrill who prompted me to think about more
general TCA's, and supplied me with several additional references.
Part of this 
work was carried out at the Mathematical Sciences Research Institute, Berkeley,
CA and hence supported by the NSF.
\vsp

\pagebreak

REFERENCES.
\vsp

[A1] F. AKMAN, A characterization of the differential in semi-infinite
cohomology, {\it J. Algebra} {\bf 162} (1993), 194-209; preprint
hep-th/9302141.
\vsp

[A2] F. AKMAN, ``The semi-infinite Weil complex of a graded Lie algebra'',
Ph.D. Thesis, Yale University, 1993.
\vsp

[A3] F. AKMAN, ``Some cohomology operators in 2-D field theory'', 
in the Proceedings of the Conference on Quantum Topology, Kansas State
University, Manhattan, KS, 24-28 March 1993, ed. David N. Yetter, World
Scientific, Singapore, 1994; preprint hep-th/9307153.
\vsp

[AKSZ] M. ALEXANDROV, M. KONTSEVITCH, A. SCHWARZ, AND O. ZABORONSKY,
The geometry of the master equation and topological quantum field theory,
preprint hep-th/9502010.
\vsp

[AD] J. ALFARO AND P.H. DAMGAARD, Non-abelian antibrackets, preprint
hep-th/9511066.
\vsp

[BV] I.A. BATALIN AND G.A. VILKOVISKY, Gauge algebra and quantization, 
{\it Phys. Lett.} {\bf 102B} (1981), 27-31.

I.A. BATALIN AND G.A. VILKOVISKY, Quantization of gauge theories with
linearly dependent generators, {\it Phys. Rev.} {\bf D28} (1983),
2567-2582.
\vsp

[BP] P. BOUWKNEGT AND K. PILCH, The BV-algebra structure of ${\cal W}_3$
cohomology, to appear in the Proceedings of ``G\"{u}rsey Memorial
Conference I: Strings and Symmetries'', eds. M. Serdaro\u{g}lu et al., 
Springer-Verlag, Berlin; preprint USC-94/17.
\vsp

[DL] C. DONG AND J. LEPOWSKY, ``Generalized vertex algebras and relative 
vertex operators'', Progress in Mathematics v. 112, Birkh\"{a}user,
Boston, 1993.
\vsp

[Fi] JOS\'{E} M. FIGUEROA-O'FARRILL, Are all TCFT's obtained by twisting
N$=2$ SCFT's?, talk given at the workshop on Strings, Gravity, and Related
Topics, held at the ICTP (Trieste, Italy) on June 29-30, 1995; preprint
hep-th/9507024.
\vsp

[FS] JOS\'{E} M. FIGUEROA-O'FARRILL AND SONIA STANCIU, Nonreductive
WZW models and their CFT's, preprint hep-th/9506151, QMW-PH-95-16.
\vsp

[FHL] I.B. FRENKEL, Y.-Z. HUANG, AND J. LEPOWSKY, On axiomatic approaches to
vertex operator algebras and modules, {\it Memoirs AMS} (1992).
\vsp

[FLM] I.B. FRENKEL, J. LEPOWSKY, AND A. MEURMAN, ``Vertex operator
algebras and the Monster'', Academic Press, New York, 1988.
\vsp

[Geb] R.W. GEBERT, Introduction to vertex algebras, Borcherds algebras,
and the Monster Lie algebra, preprint hep-th/9308151 and DESY 93-120.
\vsp

[Ger] M. GERSTENHABER, On the deformation of rings and algebras, {\it Ann.
Math.} {\bf 79} (1964), 59.

M. GERSTENHABER, The cohomology structure of an associative ring, {\it Ann.
Math.} {\bf 78} (1962), 267.
\vsp

[GV] M. GERSTENHABER AND A.A. VORONOV, Homotopy G-algebras and moduli space
operad, preprint MPI/94-71, Max-Planck-Institut in Bonn 1994, hep-th/9409063.
\vsp

[Get1] E. GETZLER, Batalin-Vilkovisky algebras and two-dimensional topological
field theories,{\it Commun. Math. Phys.} {\bf 159} (1994), 265-285; preprint
hep-th/9212043.  
\vsp

[Get2] E. GETZLER, Manin pairs and topological field theory, {\it Ann.
Phys.} {\bf 237} (1995), 161-201.
\vsp

[GHV] W. GREUB, S. HALPERIN, AND R. VANSTONE, ``Connections, curvature,
and cohomology'', v.3, Academic Press, New York, 1972-1976.
\vsp

[Gr] A. GROTHENDIECK, ``El\'{e}ments de G\'{e}ometrie Alg\'{e}brique
IV, Etude locals des sch\'{e}mas et des morphismes de sch\'{e}mas'', Pub.
Math. IHES \# 32, 1967 (Prop. 16.8.8 on p.42).
\vsp

[Hu] Y.-Z. HUANG, Operadic formulation of topological vertex algebras
and Gerstenhaber or Batalin-Vilkovisky algebras, {\em Commun. Math.
Phys.} {\bf 164} (1994), 105-144.  
\vsp

[IR] J.M. ISIDRO AND A.V. RAMALLO, Topological current algebras in
two dimensions, {\it Phys. Lett.} {\bf 316B} (1993) 488-495; preprint
hep-th/9307176, US-FT-7/93.
\vsp

[Ka] Y. KAZAMA, Novel topological field theories, {\it Mod. Phys. Lett.}
{\bf A6} (1991), 1321-1332.
\vsp

[KSV] T. KIMURA, J. STASHEFF, AND A.A. VORONOV, Homology of moduli
spaces of curves and commutative homotopy algebras, preprint alg-geom/9502006.
\vsp

[KMS] I. KOL\'{A}\v{R}, P.W. MICHOR, AND J. SLOV\'{A}K, ``Natural
Operations in Differential Geometry'', Springer-Verlag, Berlin, 1993.
\vsp

[Ko] J.-L. KOSZUL, Crochet de Schouten-Nijenhuis et cohomologie,
{\it Ast\'{e}risque} (1985), 257-271.
\vsp

[LS] T. LADA AND J.D. STASHEFF, Introduction to sh Lie algebras for
physicists, preprint hep-th/9209099, UNC-MATH-92/2.
\vsp

[Lo] J.-L. LODAY, ``Cyclic Homology'', Grundlehren der mathematischen
Wissenschaften {\bf 301}, Springer-Verlag, Berlin, 1992.
\vsp

[Li] B. H. LIAN, On the classification of simple vertex operator algebras,
{\it Commun. Math. Phys.} {\bf 163} (1994), 307-357.
\vsp

[LZ1] B. H. LIAN AND G. J. ZUCKERMAN, New perspectives on the 
BRST-algebraic structure of string theory, {\it Commun. Math. Phys.}
{\bf 154} (1993), 613-646; preprint hep-th/9211072; MR 94e:81333.
\vsp

[LZ2] B. H. LIAN AND G. J. ZUCKERMAN, Some classical and quantum algebras,
in ``Lie Theory and Geometry'', eds. Brylinski et al, {\em Progress in
Mathematics} {\bf 123}, Birkh\"{a}user, Boston, 1994, 509-529;
preprint hep-th/9404010.
\vsp

[LZ3] B. H. LIAN AND G. J. ZUCKERMAN, Commutative quantum operator
algebras, to appear in {\it J. Pure Appl. Alg.} {\bf 100} (1995); preprint
q-alg/9501014.
\vsp

[LZ4] B. H. LIAN AND G. J. ZUCKERMAN, Moonshine cohomology, preprint 
q-alg/9501015.
\vsp

[MQ] V. MATHAI AND D. QUILLEN, Superconnections, Thom classes, and equivariant
differential forms, {\it Topology} {\bf 25} (1986), 85-110.
\vsp

[Mi1] P.W. MICHOR, A generalization of Hamiltonian mechanics, {\it Journal
of Geometry and Physics} {\bf 2} (1985), 67-82.
\vsp

[Mi2] P.W. MICHOR, Remarks on the Schouten-Nijenhuis bracket, {\it
Supplemento ai Rendiconti del Circolo Matematico di Palermo, Serie II}
{\bf 16} (1987), 207-215.
\vsp

[Ni] A. NIJENHUIS, Jacobi-type identities for bilinear differential
concomitants of certain tensor fields I, {\it Koninklijke Nederlandse
Akademie van Wetenschappen, Series A, Proceedings} {\bf 58} (1955),
390-403 (journal later continued as {\it Indagationes Math.}).
\vsp

[NR] A. NIJENHUIS AND R.W. RICHARDSON, JR., Deformations of Lie algebra
structures, {\it Journal of Mathematics and Mechanics} {\bf 17} (1967),
89-105.
\vsp

[PS] M. PENKAVA AND A. SCHWARZ, On some algebraic structures arising in
string theory, preprint hep-th/9212072.
\vsp

[Sn] J.A. SCHOUTEN, \"{U}ber Differentialkonkomitanten zweier 
kontravarianter Gr\"{o}{\ss}en, {\it Koninklijke Nederlandse Akademie
van Wetenschappen, Series A, Proceedings} {\bf 2} (1940), 449-452.
Also see

J.A. SCHOUTEN, On the differential operators of first order in tensor
calculus, in Convegno Internazionale di Geometria Differenziale, Italy,
Sept. 20-26 1923, Edizioni Cremonese delle Casa Editrice Perrella,
Rome, 1954.
\vsp

[Sc] A. SCHWARZ, Geometry of Batalin-Vilkovisky quantization, preprint
hep-th/9205088.
\vsp

[St1] J. STASHEFF, Differential graded Lie algebras, quasi-Hopf algebras
and higher homotopy algebras, ...
\vsp

J. STASHEFF, Homotopy associativity of {\it H}-spaces I and II, 
{\it AMS Trans.} {\bf 108} (1963), 275-292 and 293-312.
\vsp

[St2] J. STASHEFF, Homological reduction of constrained Poisson algebras,
to appear in {\it J. Diff. Geom.}
\vsp

[Wi1] E. WITTEN, A note on the anti-bracket formalism, preprint
IASSNS-HEP-90/9.
\vsp

[Wi2] E. WITTEN, Ground ring of the two dimensional string theory,
{\em Nucl. Phys.} {\bf B373} (1992), 187.
\vsp

[Zw] B. ZWIEBACH, Closed string field theory: Quantum action and the
Batalin-Vilkovisky master equation, {\it Nucl. Phys.} {\bf B390} (1993),
33-152.
\vsp

\end{document}